\documentclass[reqno,11pt]{amsart}  

 \usepackage{gensymb}
 
 \usepackage{amsmath}
\usepackage{amssymb} 
\usepackage{mathtools}   
\usepackage{graphicx}
\usepackage{graphics}  
\usepackage{color}
\usepackage{verbatim} 
\usepackage[normalem]{ulem} 
\usepackage[mathscr]{eucal}
\usepackage{upgreek}
\usepackage{enumerate} 
 \usepackage[english]{babel} 
 \usepackage{float}
 \usepackage{subfigure}

\usepackage{babel,blindtext}

\floatstyle{boxed}
\restylefloat{figure}

\usepackage[titletoc]{appendix}
\numberwithin{equation}{section}

\usepackage[refpage,noprefix]{nomencl}
\usepackage{nomencl} 

\makenomenclature

\newtheorem{theorem}{Theorem}[section] 
\newtheorem{lemma}[theorem]{Lemma} 
\newtheorem{proposition}[theorem] {Proposition} 
 
\newtheorem{remark}[theorem]  {Remark} 
\newtheorem{definition}[theorem] {Definition}

\theoremstyle{definition}

\DeclareMathAlphabet{\mathpzc}{OT1}{pzc}{m}{it}

\newcommand{\bP} {\boldsymbol{P}}

\newcommand{\bO}{\boldsymbol{\Omega}}

\newcommand{\bbeta} {\boldsymbol{\beta}}
\newcommand{\bgamma} {\boldsymbol{\gamma}}
\newcommand{\beeta}{\boldsymbol{\eta}}
\newcommand{\bs} {\boldsymbol{\sigma}}
\newcommand{\bo}{\boldsymbol{\omega}}

\newcommand{\bzeta}{\boldsymbol{\zeta}}

\newcommand{\bDel}{\boldsymbol{\Del}}

\renewcommand{\L} {\Lambda} 
\renewcommand{\O} {\Omega}

\def\d{\delta} 
 
\newcommand{\eps}{\varepsilon}

\def\l{\lambda} 
 
\def\om{\omega}

\newcommand{\sprime}{^{\prime}\mkern-1.2mu}
\newcommand{\dprime}{^{\prime\prime}\mkern-1.2mu}

\newfam\Bbbfam 
\font\tenBbb=msbm10 
\font\sevenBbb=msbm7 
\font\fiveBbb=msbm5 
\textfont\Bbbfam=\tenBbb 
\scriptfont\Bbbfam=\sevenBbb 
\scriptscriptfont\Bbbfam=\fiveBbb

\newcommand{\R}     {\mathbb{R}} 
\newcommand{\Z}     {\mathbb{Z}} 
\newcommand{\N}     {\mathbb{N}} 
\renewcommand{\P}   {\mathbb{P}} 
 
\newcommand{\E}     {\mathbb{E}}

 \newcommand{\floor}[1]{\left\lfloor #1 \right\rfloor}

\def\1{{\mathchoice {1\mskip-4mu\mathrm l}      
{1\mskip-4mu\mathrm l} 
{1\mskip-4.5mu\mathrm l} {1\mskip-5mu\mathrm l}}} 
\newcommand{\ssup}[1] {{\scriptscriptstyle{({#1}})}} 
\def\comment#1{} 
\newtheoremstyle{thm}{2ex}{2ex}{\itshape\rmfamily}{} 
{\bfseries\rmfamily}{}{1.7ex}{} 
 
\newtheoremstyle{rem}{1.3ex}{1.3ex}{\rmfamily}{} 
{\itshape\rmfamily}{}{1.5ex}{} 
 
\newenvironment{proofsect}[1] 
{\vskip0.1cm\noindent{\scshape #1.}\hskip0.5cm}

\newcommand{\Fcal}   {{\mathcal F }}

\newcommand{\Mcal}   {{\mathcal M }}

\newcommand{\Ucal}   {{\mathcal U }}

\newcommand{\Lscr} {\mathscr{L}}
\newcommand{\Zscr}{\mathscr{Z}}

\newcommand{\Escr}{{\mathscr{E}}}

 \newcommand{\ex}{{\rm e}} 
 
\renewcommand{\d}{{\rm d}}

\newcommand{\Leb}{{\rm Leb}}

\newcommand{\dist}{{\operatorname {dist}}} 
\newcommand{\area}{{\operatorname {\mathbf{area}}}}

\newcommand{\radius}{{\operatorname {radius}}} 
 
\newcommand{\Del}{{\operatorname{\sf Del}}}
\newcommand{\Vor}{{\operatorname {\sf Vor}}}
\newcommand{\argmax}{{\operatorname {argmax }}}

\newcommand{\Exp}{\mathscr{E}\kern-0.2mm{\operatorname{xp}}}
\newcommand{\Log}{\mathscr{L}\kern-0.2mm{\operatorname{og}}}

\newcommand{\heap}[2]{\genfrac{}{}{0pt}{}{#1}{#2}} 
\newcommand{\pr}{{\operatorname {pr}}}
\renewcommand{\emptyset} {\varnothing}

\makeindex

\begin{document}

\title[\hfill Phase transitions\hfill]
{The Widom-Rowlinson Model on the Delaunay Graph}


\author{Stefan Adams and  Michael Eyers}
\address{Mathematics Institute, University of Warwick, Coventry CV4 7AL, United Kingdom}
\email{S.Adams@warwick.ac.uk}

\thanks{}
  

 


\begin{abstract}

We establish phase transitions for continuum Delaunay multi-type particle systems (continuum Potts or Widom-Rowlinson models) with  a repulsive interaction between particles of different types. Our interaction potential depends solely on the length of the Delaunay edges. We show that a phase transition occurs for sufficiently large  activities and for sufficiently large potential parameter   proving  an old conjecture of Lebowitz and Lieb extended to the Delaunay structure. 
Our approach involves a Delaunay random-cluster representation analogous to the Fortuin-Kasteleyn representation of the Potts model. The phase transition manifests itself  in the mixed site-bond percolation of the corresponding random-cluster model. Our proofs rely mainly on geometric properties of  Delaunay tessellations in  $\R^2 $ and on recent studies \cite{DDG12} of Gibbs measures for geometry-dependent interactions. The main tool is a uniform bound on the number of connected components in the Delaunay graph which provides a novel approach to Delaunay Widom Rowlinson models based on purely geometric arguments. The interaction potential ensures that shorter Delaunay edges are more likely to be open and thus offsets the possibility of having an unbounded number of connected components. 
\end{abstract} 
\maketitle


\section{Introduction and results}

\subsection{Introduction}
Although the study of phase transitions is one of the main subjects of mathematical statistical mechanics, examples of models exhibiting phase transition are mainly restricted  to lattice systems. In the continuous setting results are much harder to obtain, e.g., the proof of a liquid-vapor phase transition in \cite{LMP99}, or the spontaneous breaking of rotational symmetry in two dimensions for a Delaunay hard-equilaterality like interaction \cite{MR09}. These phase transitions manifest themselves in  breaking of a continuous symmetry. There is another  specific model for which a phase transition is known to occur: the model of Widom and Rowlinson \cite{WR70}. This is a multi-type particle system in $ \R^d, d\ge 2 $, with hard-core exclusion between particles of different type, and no interaction between particles of the same type. The phase transition in this model was stablished by Ruelle \cite{Rue71}. Lebowitz and Lieb \cite{LL72} extended his result by replacing the hard-core exclusion by a soft-core repulsion between unlike particles. Finally, phase transition results for classes of continuum Potts models in $ \R^d, d\ge 2 $, have been derived in \cite{GH96}. The phase transitions for large activities in all these systems reveal themselves in breaking of the symmetry in the  type-distribution.  In \cite{BBD04} the soft repulsion in \cite{GH96} between unlike particles has been replaced by another kind of soft repulsion based on the structure of some graph. 

In this paper we establish the existence of a phase transition  for a class of continuum Delaunay Widom-Rowlinson (Potts) models in $ \R^2 $. The repulsive interaction between unlike particles is of finite range, and it depends  on the geometry of the Delaunay tessellation, i.e., the length of the edges.  The potential is formally given as
$$
\phi_\beta(\ell)=\log\Big(\frac{\ell^4+\beta}{\ell^4}\Big)\1\{\ell\le R\},\quad \beta>0, \ell\ge 0,
$$ where $ \ell\ge 0 $ is the length of an Delaunay edge and $ \beta>0 $ is the potential parameter and $ R>0 $ is finite range condition of the potential. The main novelty of our paper is a uniform bound on the number of connected components in the Delaunay random cluster model which is purely based on geometrical properties of Delaunay tessellations in two dimensions.  The potential ensures that for large values of the parameter $ \beta >0 $,  Delaunay edges with shorter lengths are more likely to be connected than those with longer edges, enabling us to bound the number of connected components for clouds of points with vanishing point-wise distances.   This paper is an extensive further development of the recent work \cite{AE16} where all models had an additional background hard-core potential introducing a length scale for the configurations. Gibbs models on Delaunay structures have been studied in \cite{BBD99,BBD02,BBD04,Der08,DDG12,DG09,DL11}, and our results rely on the existence of Gibbs measures for the geometry-dependent interaction  using methods in \cite{DDG12}. Our approach is based on a Delaunay random-cluster representation.   A phase transition for our Delaunay Potts model follows if we can show that the corresponding percolation process contains an infinite cluster. A similar program was carried out by Chayes et al. in \cite{CCK} for the hard-core Widom-Rowlinson model. In that case, the existence of infinite clusters follows from a stochastic comparison with the Poisson Boolean model of continuum percolation, while our framework uses a coarse graining method to derive a stochastic comparison with mixed site-bond  percolation on $ \Z^2 $. Our results are extension of \cite{LL72} and \cite{CCK} to the Delaunay structure replacing hard-core constraint by our soft-core repulsion. In particular we obtain phase transition for all activities once the interaction parameter $ \beta>0 $ (inverse temperature) is sufficiently  large depending on the activity. We note that our random-cluster representation requires the symmetry of the type  interaction. In the non-symmetric Widom-Rowlinson models, the existence of a phase transition has been established by Bricmont et al. \cite{BKL}, and recently by Suhov et al. \cite{MSS}.

\subsection{Remarks on Delaunay tessellations}
We add  some remarks on models defined on Delaunay hypergraph structures.  There are differences between geometric models on the Delaunay hypergraph structure and  classical particle models such as the Widom-Rowlinson model \cite{WR70}  and its soft-core variant of Lebowitz and Lieb \cite{LL72}. The first is that edges and triangles in the Delaunay hypergraph are each proportional in number to the number of particles in the configuration. However, in the case of the complete hypergraph the number of edges is proportional to the number of particles squared and the number of triangles is proportional to the number of particles cubed. Secondly, in complete graphs of all classical models, the neighbourhood of a given point depends only on the distance between points and so the number of neighbours increases with the activity parameter $ z $ of the underlying point process. This means that the system will become strongly connected for high values of $z$. This is not the case for the Delaunay hypergraphs which exhibit a self-similar property. Essentially, as the activity parameter $z$ increases, the expected number of neighbours to a given point in the Delaunay hypergraph remains the same, see \cite{M94}. Therefore, in order to keep a strong connectivity, we use a type interaction between particles of Delaunay edges  with a non-constant mark. Finally, and perhaps most importantly, is the question of additivity. Namely, suppose we have an existing particle configuration $\om$ and we want to add a new particle $x$ to it. In the case of classical many-body interactions, this addition will introduce new interactions that occur between $x$ and the existing configuration $\om$. However, the interactions between particles of $\om $ remain unaffected, and so classical many-body interactions are additive. On the other hand, in the Delaunay framework, the introduction of a new particle to an existing configuration not only creates new edges and triangles, but destroys some too. The Delaunay interactions are therefore not additive, and for this reason, attractive and repulsive interactions are indistinct. In the case of a hard exclusion interaction, we arrive at the possibility that a configuration $ \om $ is excluded, but for some $x$, $ \om\cup x $ is not. This is called the non-hereditary property \cite{DG09}, which seems to rule out using techniques such as stochastic comparisons of point processes \cite{GK97}.

\subsection{Setup}
We consider configurations of points in $ \R^2 $ with internal degrees of freedom, or marks. Let $ M_q=\{1,\ldots,q\}, q\in\N, q\ge 2 $, be the  finite set of different marks. That is, each marked point is represented by a position $ x\in\R^2 $ and a mark $ \sigma(x)\in M_q $, and each marked configuration $ \bo $ is a countable subset of $ \R^2\times M_q $ having a locally finite projection onto $ \R^2 $. We denote by $ \bO $ the set of all marked configurations with locally finite projection onto $ \R^2 $. We will sometimes identify $ \bo $ with a vector $ \bo=(\omega^{\ssup{1}},\ldots,\om^{\ssup{q}}) $ of pairwise disjoint locally finite sets $ \om^{\ssup{1}},\ldots,\om^{\ssup{q}} $ in $ \R^2 $ (we write $ \O $ for the set of all locally finite configurations in $ \R^2 $). Any $ \bo $ is uniquely determined by the pair $ (\om,\sigma)$, where $ \omega=\cup_{i=1}^q\omega^{\ssup{i}} $ is the set of all occupied positions, and where the mark   function $ \sigma\colon \om\to M_q $ is defined by $ \sigma(x)=i $ if $ x\in\om^{\ssup{i}}, i\in M_q $. For each measurable set $B$  in $ \R^2\times M_q $ the counting variable $ N(B)\colon\bo\to\bo(B) $ on $ \bO $ gives the number of marked particles such that  the pair (position, mark) belongs to $B$. We equip the space $ \bO $ with the $\sigma$-algebra $\boldsymbol{\Fcal} $ generated by the counting variables $N(B) $ and the space  $ \O $ of locally finite configurations with the $\sigma$-algebra $ \Fcal $ generated by the counting variables $ N_\Delta=\#\{\om\cap\Delta\} $ for $ \Delta\Subset\R^2 $ where we write $ \Delta\Subset\R^2 $ for any bounded $ \Delta\subset\R^2 $. As usual, we take as reference measure on $ (\bO,\boldsymbol{\Fcal}) $ the marked Poisson point process $ \boldsymbol{\Pi}^z $ with intensity measure $ z\Leb\otimes\mu_{\sf u} $  where $ z>0 $ is an arbitrary activity, $ \Leb $ is the Lebesgue measure in $ \R^2 $, and $ \mu_{\sf u} $ is the uniform probability measure on $ M_q  $.

For each $ \L\subset\R^2 $ we write $ \bO_\L=\{\bo\in\bO\colon \bo\subset\L\times M_q\} $  for the set of configurations in $ \L $, $  \pr_\L\colon\bo\to\bo_\L:=\bo\cap \L\times M_q $ for the projection from $ \bO $ to $ \bO_\L $ (similarly for unmarked configurations), $\boldsymbol{\Fcal}_\L^\prime=\boldsymbol{\Fcal}|_{\bO_\L} $ for the trace $\sigma$-algebra of $ \boldsymbol{\Fcal} $ on $ \bO_\L $, and  $ \boldsymbol{\Fcal}_\L=\pr_\L^{-1}\boldsymbol{\Fcal}_\L^\prime\subset\boldsymbol{\Fcal} $ for the $\sigma$-algebra of all events that happen in $ \L $ only. The reference measure on $ (\bO_\L,\boldsymbol{\Fcal}_\L^\prime) $ is $ \boldsymbol{\Pi}_\L^z:=\boldsymbol{\Pi}^z\circ\pr_\L^{-1} $. In a similar way we define the corresponding objects for unmarked configurations, $ \Pi^z, \Pi^z_\L, \O_\L, \pr_\L, \Fcal_\L^\prime $, and $ \Fcal_\L $. Finally, let $ \varTheta=(\vartheta_x)_{x\in\R^2} $ be the shift group, where $ \vartheta_x\colon\bO\to\bO $ is the translation of the  spatial component by the vector $ -x\in\R^2 $. Note that by definition, $ N_\Delta(\vartheta_x\bo)=N_{\Delta+x}(\bo) $ for all $ \Delta\subset\R^2 $. \medskip

We outline the definitions for the unmarked configurations first with obvious adaptations to the case of marked point configurations.  The set $ \Del $ of Delaunay hyperedges consist of all pairs $ (\eta,\om)$ with $ \eta\subset\om $ for which there exits an open ball $ B(\eta,\om) $ with $ \partial B(\eta,\om)\cap\om=\eta $ that contains no points of $ \om $. For $ m=1,2,3  $, we write $\Del_m=\{(\eta,\om)\in\Del\colon \#\eta=m\} $ for the set of Delaunay simplices with $ m $ vertices. Given a configuration $\om $  the set of all Delaunay hyperedges $ \eta\subset\om $ with $ \#\eta=m $ is denoted by $\Del_m(\om) $. It is possible that $ \eta\in\Del(\om) $ consists of four or more points on a sphere with no points inside. In fact, for this not to happen, we must consider configurations in general   position as in \cite{M94}. More precisely, this means that no four points lie on the boundary of a circle and every half-plane contains at least one point. Fortunately, this occurs with probability one for our Poisson reference measure, and in fact, for any stationary point process. Note that the open ball $ B(\eta,\om) $ is only uniquely determined when $ \#\eta=3 $ and $ \eta $ is affinely independent. Henceforth, for each configuration $ \om $ we have an associated Delaunay triangulation 
\begin{equation}\label{triangulation}
\{\tau\subset \om\colon \#\tau=3,B(\tau,\om)\cap\om=\emptyset\}
\end{equation}
of the plane, where $ B(\tau,\om) $ is the unique open ball with $ \tau\subset\partial B(\tau,\om) $. The set in \eqref{triangulation} is uniquely determined and defines a triangulation of the convex hull of $ \om $ whenever $ \om $ is in general position (\cite{M94}). In a similar way one can define the marked Delaunay hyperedges, $ \boldsymbol{\Del} $ and $ \boldsymbol{\Del}_m(\bo) $ respectively,   where the Delaunay property refers to the spatial component only.

\medskip

Given a configuration $ \om \in\Omega $ (or $ \bo$) we write $ \O_{\L,\om}=\{\zeta\in\O\colon \zeta\setminus\L=\om\} $ (resp. $ \bO_{\L,\bo} $) for the set of configurations which equal $ \om $ off $ \L $.  For any edge $ \eta\in\Del_2 $ we denote its  length by $ \ell(\eta) $, i.e., $ \ell(\eta)=|x-y| $ if $ \eta=\{x,y\} $. The interaction is given by the following Hamiltonian in $ \L $ with boundary condition $ \bo\in\bO $,  
\begin{equation}\label{Hamiltonian}
H_{\L,\bo}(\bzeta):=\sum_{\heap{\beeta\in\bDel_{2,\L}(\bzeta)\colon}{\eta\in\Del_{2,\L}(\zeta)}}\phi_\beta(\ell(\eta))(1-\delta_\sigma(\beeta)),\quad\bzeta\in\O_{\L,\bo},
\end{equation}
where $ \Del_{2,\L}(\zeta):=\{\eta\in\Del_2(\zeta)\colon \,\exists\,  \tau\in\Del_3(\zeta), \eta\subset\tau,\partial B(\tau,\zeta)\cap\L\not=\emptyset\} $.
Here   $ \phi_\beta $ is a measurable function of the length $ \ell(\eta) $ of an edge defined for any $\beta\ge 0$,
\begin{equation}
\phi_\beta(\ell)=\log\Big(\frac{\ell^4+\beta}{\ell^4}\Big)\1\{\ell\le R\},
\end{equation}
and 
$$ 
\delta_\sigma(\beeta)=\begin{cases} 1 &, \mbox{ if } \sigma_{\beeta}(x)=\sigma_{\beeta}(y) \mbox{ for } \eta=\{x,y\},\\
0&, \mbox{ otherwise}.
\end{cases}
$$ 

\noindent Note the following scaling relation for the potential
\begin{equation}
\phi_\beta(L\ell)=\phi_{\beta/L^4}(\ell),\qquad\mbox{for any } L>0 \;\mbox{ with } L\ell\le R, \ell\le R.
\end{equation}
\smallskip

 Following \cite{DDG12} we say a  configuration $ \bo\in\bO $   (or $\om\in\O$) is \textbf{admissible} for $ \L\Subset\R^2 $ and activity $ z$ if $ H_{\L,\bo} $ is $\boldsymbol{\Pi}^z$-almost surely well-defined and $ 0<Z_\L(\bo)<\infty $, where the  partition function is defined as
  $$
  Z_\L(\bo)=\int_{\bO_{\L,\bo}}\,{\rm e}^{-H_{\L,\bo}(\bzeta)}\,\Pi^z(\d\zeta^{\ssup{1}})\cdots\Pi^z(\d\zeta^{\ssup{q}}).
  $$
 We denote the set of admissible configurations by $ \bO_\L^* $. The Gibbs distribution for $ \phi_\beta $, and $ z> 0 $ in $\L $ with admissible boundary condition $ \bo $ is defined as
 \begin{equation}\label{Gibbsdist}
 \gamma_{\L,\bo}(A)=\frac{1}{Z_\L(\bo)}\int_{\bO_{\L,\bo}}\,\1_A(\bzeta\cup\bo){\rm e}^{-H_{\L,\bo}(\bzeta)}\,\Pi_\L^z(\d\bzeta),\quad A\in\boldsymbol{\Fcal}.
 \end{equation}
  It is evident from \eqref{Gibbsdist} that, for fixed $ \zeta\in\O_\L $, the conditional distribution of the marks of $ \bzeta= (\zeta^{\ssup{1}},\ldots,\zeta^{\ssup{q}})$  relative to $ \gamma_{\L,\bo} $ is that of a discrete Potts model on $ \zeta $ embedded in the Delaunay triangulation with position-dependent interaction between the marks. This justifies calling our model  \textit{Delaunay Potts model or Delaunay Widom-Rowlinson model}.
  
  \begin{definition}
  A probability measure $ \mu $ on $ \bO $ is called a Gibbs measure for the Delaunay  Potts model with activity $ z>0 $ and interaction potential $ \phi_\beta $  if $ \mu(\bO^*_\L)=1 $ and
 \begin{equation}\label{DLR}
   \E_\mu[f]=\int_{\bO^*_\L}\,\frac{1}{Z_\L(\bo)}\int_{\bO_{\L,\bo}}f(\bzeta\cup\bo){\rm e}^{-H_{\L,\bo}(\bzeta)}\,\boldsymbol{\Pi}^z_\L(\d\bzeta)\mu(\d\bo)
 \end{equation} for every $ \L\Subset\R^2 $ and every measurable function $ f$.
  \end{definition}
  The equations in \eqref{DLR} are the DLR equations (after Dobrushin, Lanford, and Ruelle). They ensure that the Gibbs distribution in \eqref{Gibbsdist} is a version of the conditional probability $ \mu(A|\boldsymbol{\Fcal}_{\L^{\rm c}})(\bo) $. The measurability of all objects is established in \cite{E14,DDG12}.

\subsection{Results and remarks}

  \begin{proposition}[\textbf{Existence of Gibbs measures}]\label{Propexist}
  
 For any $z>0 $ there exist at least one Gibbs measure for the Delaunay Widom-Rowlinson (Potts)  model     with parameter $ \beta >  0$.  
 
  \end{proposition}

 \begin{remark}[\textbf{Gibbs measures}]
 The proof is using the so-called pseudo-periodic configurations (see Appendix~\ref{pseudoperiodic} or \cite{DDG12}) and properties of the potential $ \phi_\beta $. Existence of Gibbs measures for related different  Delaunay models have been obtained in \cite{BBD99, Der08,DG09}. Note that for $q=1$ our models have no marks and Gibbs measures do exist as well (\cite{DDG12}).\hfill $ \diamond$
 \end{remark}
 
A \textit{phase transition}  is said to occur if there exists more than one Gibbs measure for the Delaunay Potts model. The following  theorem shows that this happens for all activities $ z $ and sufficiently large parameter $ \beta $ depending on $z$.  Note that $ \beta $ is a parameter for the type interaction and not the usual inverse temperature.

  \begin{theorem}[\textbf{Phase transition}]\label{THM-mainedge}
  For all $ \ell\in(0,\frac{R}{2\sqrt{3}}] $ and $ \beta> q $ there exist $ \alpha^*=\alpha^*(R,q) $ and  $ z_0=z_0(\alpha^*,q,\ell) \ge z_0^*(\alpha^*,q) $ such that for all $ z\ge z_0 $ there exists $ \beta_0=\beta_0(q,R,z) $ such that for all $ \beta\ge\beta_0\vee q $  there exit at least $q$ different  Gibbs measures for the Delaunay Widom Rowlinson (Potts) model. 
  \end{theorem}

\begin{remark}
Theorem~\ref{THM-mainedge} also holds for any potential depending on the length $ \ell $ of Delaunay edges
$$ 
\phi_\beta^{\ssup{\gamma}}(\ell):=\log\Big(\frac{\ell^{3+\gamma}+\beta}{\ell^{3+\gamma}}\Big)\1\{\ell\le R\},\qquad \gamma>0.
$$ \hfill $ \diamond $
\end{remark}

\begin{remark}[\textbf{Free energy}]
One may wonder if the phase transition manifest itself thermodynamically by a non-differentiability ("discontinuity") of the free energy (pressure).  Using the techniques from \cite{Geo94} and  \cite{DG09}, it should be possible to obtain a variational representation of the free energy, see also \cite{ACK11} for free energy representations for marked configurations. Then a discontinuity of the free energy can be established using our results above. For continuum Potts models this has been established in \cite[Remark~4.3]{GH96}.    \hfill $ \diamond$
\end{remark}

 \begin{remark}[\textbf{Uniqueness of Gibbs measures}]
 To establish uniqueness of the Gibbs measure in our Delaunay Potts model one can use the Delaunay random-cluster measure $ C_{\L_n,\omega} $, to be defined in \eqref{DRCmeasure} below. In \cite[Theorem~6.10]{GHM} uniqueness is established once the probability of an open connection of the origin to infinity is vanishing for the limiting lattice version of the random-cluster measure, that is, for some set $ \Delta\Subset\R^2$ containing the origin, 
 $$
 \lim_{n\to\infty} C_{\L_n,\om}(\Delta\longleftrightarrow \L_n^{\rm c})=0,
 $$  for a sequence of boxes $ \L_n\Subset \R^2$ with $ \L_n\uparrow\R^2 $ as $ n\to\infty $. One way to achieve this, is to obtain an stochastic domination of the Delaunay  random-cluster measure by the so-called random Delaunay edge model of hard-core particles.  Using \cite{BBD02} we know that the critical probabilities for both, the site and bond percolation on the Delaunay graph, are bounded from below. Extension to our Delaunay edge percolation can provide a corresponding lower bound as well. Thus, if our parameter $ \beta $ is chosen sufficiently small, then there is no percolation in our Delaunay random-cluster measure and therefore uniqueness of the Gibbs measure.  \hfill $ \diamond$ 
 \end{remark}  
  
 The study for Widom-Rowlinson or Potts models with geometry-dependent interaction is by far not complete, one may wish to extend the single edge  (or triangle) interaction to mutual adjacent Voronoi cell interaction.  The common feature of all these ``ferromagnetic'' systems is that phase transitions are due to breaking the symmetry  of the type  distribution.  
 
 The rest of the paper is organised a follows. In Section~\ref{DelCluster} we define the Delaunay random-cluster measure for edge configruations, and in Section~\ref{tileperc} we establish percolation in this model for certain parameters. The main novelty is the extensive and elaborate proof of the uniform bound on the number of connected components using purely geometric properties in Section~\ref{number}. Finally, in Section~\ref{proofs} we gives details of our remaining proofs.  
 
 \section{The random cluster method}
In Section~\ref{DelCluster} we introduce the Delaunay random cluster model and show percolation for this model in Section~\ref{tileperc} via comparison with mixed site-bond percolation on $ \Z^2 $. We conclude in Section~\ref{phase} with our proof of Theorem~\ref{THM-mainedge}. 
The key step is our novel uniform estimate of the number of connected components  in Section~\ref{number}.  The proof of this bound uses solely geometric arguments and constitutes a major part of this work.

\subsection{Delaunay Random Cluster measure}\label{DelCluster}
 For $ \L\Subset\R^2 $  and parameters $ z $ and $ \phi_\beta  $ we define a joint distribution of the Delaunay Potts model and an  edge process which we call Delaunay random-cluster model. The basic idea is to introduce random edges  between points in the plane.   Let 
 $$
 E_{\R^2}=\{\eta=\{x,y\}\subset\R^2\colon x\not=y\}
 $$  be the set of all possible edges  of points in $ \R^2 $, likewise, let $E_\L $ be  the set of all edges in $ \L $ and $ E_\zeta $ for the  set of edges in $ \zeta\in\Omega_{\L,\om} $. We identify $ \om $ with $\om^{\ssup{1}} $ and $ \bo=(\om^{\ssup{1}},\emptyset,\ldots,\emptyset) $.  This allows only monochromatic boundary conditions whereas the general version involves the so-called Edwards-Sokal coupling (see \cite{GHM} for lattice Potts models). We restrict ourself to the former case for ease of notation. We write $$ \Escr =\{E\subset E_{\R^2}\colon E\mbox{ locally finite}\}$$ for the set of all locally finite edge  configurations.

 The joint distribution is built from the following three components.

 \noindent The \textit{point distribution}  is  given by the Poisson process  $\Pi^{zq} $ for any admissible boundary condition $ \om\in\O_{\L}^*$ and activity $ zq$.
  
 \medskip
 
 \noindent The \textit{type picking mechanism} for a given configuration $ \zeta\in\Omega_{\L,\om} $ is the distribution $ \l_{\zeta,\L} $ of the mark vector $\sigma\in M_q^\zeta $.  Here $ (\sigma(x))_{x\in\zeta} $ are independent and uniformly distributed random variables on $M_q$ with $ \sigma(x)=1 $ for all $ x\in\zeta_{\L^{\rm c}}=\om $.  The latter condition ensures that all points outside of $ \L$ carry the given fixed mark.
 \medskip
 
 \noindent The \textit{edge drawing mechanism}. Given a point configuration $ \zeta\in\Omega_{\L,\om} $,   we let $ \mu_{\zeta,\L} $ be the distribution of the random edge configuration $ \{\eta\in E_\zeta\colon\upsilon(\eta)=1\} \in\Escr $ with the edge configuration $ \upsilon\in \{0,1\}^{E_\zeta} $ having
 probability   
 $$
\prod_{\eta\in E_\zeta}p(\eta)^{\upsilon(\eta)}(1-p(\eta))^{1-\upsilon(\eta)}
$$ with
\begin{equation}\label{tiledrawing}
p(\eta):=\P(\upsilon(\eta)=1)=\begin{cases} (1-\ex^{-\phi_\beta(\eta)})\1_{\Del_2(\zeta)}(\eta) & \mbox{ if } \eta\in E_{\R^2}\setminus E_{\L^{\rm c}},\\
\1_{\Del_2(\zeta)}(\eta) & \mbox{ if } \eta\in E_{\L^{\rm c}}.
\end{cases} 
\end{equation}
 
The measure $ \mu_{\zeta,\L} $ is a point process on $E_{\R^2} $. Note that $ \zeta\rightarrow\l_{\zeta,\L} $ and $ \zeta\rightarrow\mu_{\zeta,\L} $ are probability kernels (see \cite{E14,AE16}). Let the measure
 $$
 P^{zq}_{\L,\bo}(\d\bzeta,\d E)=\frac{1}{Z_\L(\om)}P_{\L,\om}^{zq}(\d\zeta)\l_{\zeta,\L}(\d\bzeta)\mu_{\zeta,\L}(\d E)
 $$ be supported on the set of all $ (\bzeta,E) $ with $ \bzeta\in\bO_{\L,\bo} $ and $ E\subset E_\zeta $. We shall condition on the event that the marks of the points are constant on each connected component in the graph $ (\zeta,E\cap E_\zeta) $. Two distinct vertices $x $ and $y$ are adjacent to one another if there exists $ \eta\in E_\zeta $ such that $ \{x,y\}=\eta $. A path in the graph $ (\zeta,E\cap E_\zeta) $ is an alternating sequence $ v_1,e_1,v_2,e_2,\ldots $ of distinct vertices $ v_i $ and edges $e_j $ such that $ \{v_i,v_{i+1}\}=e_i $ for all $ i\ge 1 $.  We write
 $$
 A=\{(\bzeta,E)\in\bO\times\Escr\colon\sum_{\eta\in E}(1-\delta_{\sigma}(\beeta))=0\}
 $$ for the set of marked  point configurations such that all vertices of the edges carry the same mark. The set $ A $ is measurable which one can see from writing the condition in the following way
 $$
 \sum_{\eta=\{x,y\}\in E}\;\sum_{i=1}^q\big(\1_{\zeta^{\ssup{i}}}(x)(1-\1_{\zeta^{\ssup{i}}}(y))\big)=0
 $$ and using the fact that $ (\bzeta,x)\mapsto \1_{\zeta^{\ssup{i}}}(x) , i=1,\ldots, q $, are measurable (see \cite[Chapter 2]{GH96}). Furthermore,  $ \Pi^{zq}_{\L,\bo}(A)>0 $, which follows easily observing
 $ \Pi^{zq}_{\L,\bo}(A)\ge \Pi^{zq}_{\L,\om}(\{\widetilde{\om}\})=\ex^{-zq|\L|}/Z_\L(\om) $,  where $ \widetilde{\om} $ is the configuration  which  equals $ \om $ outside of $ \L $ and which is empty inside $ \L $.  Henceforth, the random-cluster measure
$$
 \boldsymbol{P}=P^{zq}_{\L,\bo}(\cdot\vert A)
$$
 is well-defined. We obtain the following two measures from the  random-cluster measure $\bP$, namely if we disregard the edges we obtain the Delaunay Gibbs distribution $ \gamma_{\L,\bo} $ in \eqref{Gibbsdist} (see \cite{E14}). For the second measure consider the  mapping $ \mbox{sp}\colon(\bzeta,E)\to(\zeta,E) $ from $ \bO\times\Escr $ onto $ \O\times\Escr $ where $ \bzeta\mapsto \zeta=\cup_{i=1}^q\zeta^{\ssup{i}} $. For each $ (\zeta,E) $ with $ E\subset E_\zeta $ we let $ K(\zeta,E)$ denote the number of connected components in the graph $ (\zeta,E) $. The Delaunay random-cluster distribution on $ \O\times\Escr $ is defined by
\begin{equation}\label{DRCmeasure}
  C_{\L,\om}(\d\zeta,\d E)=\frac{1}{\Zscr_\L(\om)}q^{K(\zeta,E)}\Pi^z_{\L,\om}(\d\zeta)\mu_{\zeta,\L}(\d E),
 \end{equation}
 where $ \Pi^z_{\L,\om} $ is the Poisson process with activity $ z$ replacing  $zq$ and where
 $$
 \Zscr_\L(\om)=\int_{\O_{\L,\om}}\,\int_{\Escr}\, q^{K(\zeta,E)}\Pi_{\L,\om}^z(\d\zeta)\mu_{\zeta,\L}(\d E)
 $$ is the normalisation. It is straightforward to  show that $ \bP\circ\mbox{sp}^{-1}=C_{\L,\om} $.  
 
\begin{figure}[t]
\centering
\subfigure[$ \Del_2(\zeta)$]{
\includegraphics[scale=0.8]{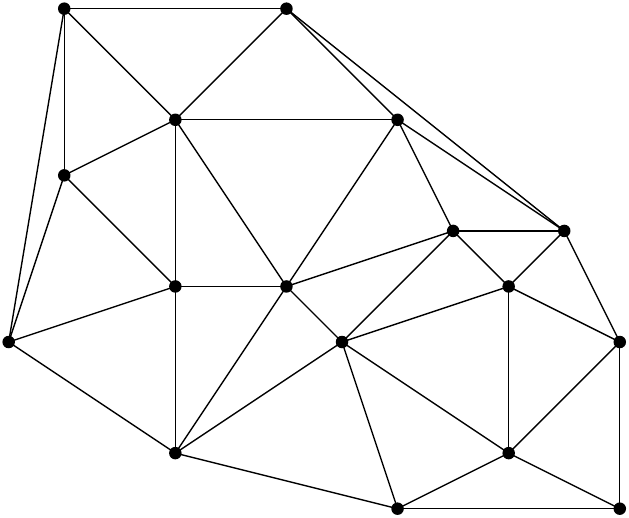}}
\subfigure[$\Del_2(\zeta\cup\{x_0\}) $]{
\includegraphics[scale=0.8]{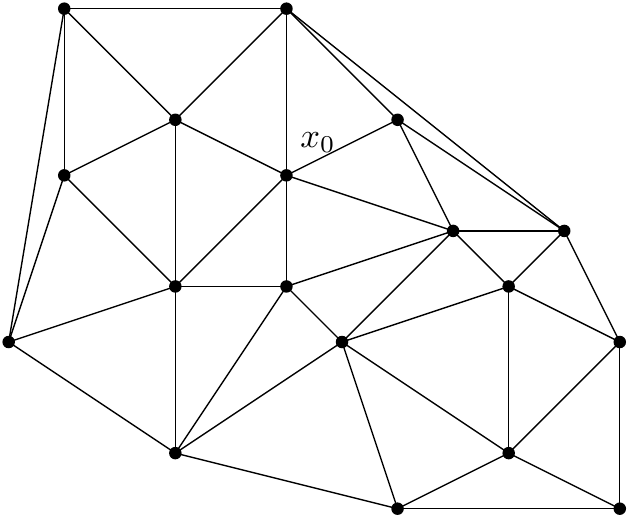}
}

\bigskip
\bigskip

\subfigure[$ E_{x_0,\zeta}^{\ssup{\rm ext}} $ ]{
\includegraphics[scale=0.8]{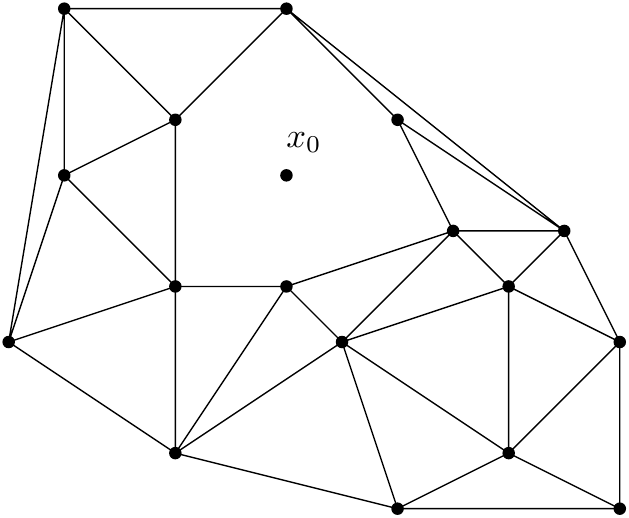}
}

\caption{The Delaunay sets $ \Del_2(\zeta), \Del_2(\zeta\cup\{x_0\}) $, and the exterior Delaunay set $ E_{x_0,\zeta}^{\ssup{\rm ext}}$ }
\label{fig}
\end{figure}

For our main proofs  we need to investigate the geometry of the Delaunay tessellation, and in particular what happens when we augment $ \zeta \in\Omega_{\L,\omega} $ with a new point $ x_0\notin\zeta $. Some edges (triangles)  may be destroyed, some are created, and some remain. This process is well described in \cite{Lis94}. We give a brief account here for the convenience of the reader. We insert the point $ x_0 $ into one of the triangles $ \tau $ in $ \Del_3(\zeta) $. We then create three new edges that join $ x_0 $ to each of the three vertices of $ \tau $. This creates three new triangles, and destroys one. We now need to verify that the new triangles each satisfy the Delaunay condition \eqref{triangulation}, that is, that their circumscribing balls contain no points of $ \zeta$. If this condition is satisfied the new triangle remains, if it is not satisfied, then there is a point $ x_1\in\zeta $  inside the circumscribing ball. We remove the edge not connected to $ x_0 $, and replace it by an edge connecting $x_0 $ and $ x_1 $. This results in the creation of two new triangles. Each of these triangles must be checked as above and the process continues. Once all triangles satisfy the Delaunay condition, we arrive at the Delaunay triangulation $ \Del_3(\zeta\cup\{x_0\}) $ and their Delaunay edges  $ \Del_2(\zeta\cup\{x_0\}) $. Let
\begin{equation}\label{tilesets}
\begin{aligned}
E_{x_0,\zeta}^{\ssup{\rm ext}}&=\Del_2(\zeta)\cap\Del_2(\zeta\cup\{x_0\}),\\
E^{\ssup{+}}_{x_0,\zeta}&=\Del_2(\zeta\cup\{x_0\})\setminus\Del_2(\zeta)=\Del_2(\zeta\cup\{x_0\})\setminus E_{x_0,\zeta}^{\ssup{\rm ext}},\\
E^{\ssup{-}}_{x_0,\zeta}&=\Del_2(\zeta)\setminus\Del_2(\zeta\cup\{x_0\})=\Del_2(\zeta)\setminus E_{x_0,\zeta}^{\ssup{\rm ext}},
\end{aligned}
\end{equation}
be the set of exterior, created, and destroyed Delaunay edges  respectively, see Figure~\ref{fig}. 

%

Note that any new triangle must contain $ x_0 $, i.e.,
$$
E^{\ssup{+}}_{x_0,\zeta}=\{\tau\in\Del_3(\zeta\cup\{x_0\})\colon \tau\cap x_0=x_0\}.$$
We let $ \mu^{\ssup{-}}_{x_0,\zeta},\mu^{\ssup{+}}_{x_0,\zeta} $, and $ \mu^{\ssup{\rm ext}}_{x_0,\zeta} $ be the edge drawing mechanisms on $ E_{x_0,\zeta}^{\ssup{\rm ext}}, E^{\ssup{+}}_{x_0,\zeta} $, and $ E^{\ssup{-}}_{x_0,\zeta}$,  respectively, which are derived from the edge drawing measure $ \mu_{\zeta,\L} $ above. The crucial step is an  estimate on the number of connected components in a neighbourhood of the point $ x_0 $.  For any $ \zeta\in\Omega_{\L,\omega} $ the neighbourhood  of a point $ x_0\in\L,x_0\notin\zeta, $ is the following random graph $ G_{x_0,\zeta}=(V_{x_0,\zeta},E_{x_0,\zeta}^{\ssup{\rm nbd}}) $ where $ V_{x_0,\zeta} $ is the set of points that share an edge with $ x_0 $ in $ \Del_2(\omega\cup\{x_0\}) $ and $ E_{x_0,\zeta}^{\ssup{\rm nbd}} $ is the set of edges in $ \Del_2(\omega\cup\{x_0\}) $ that have both endpoints in $ V_{x_0,\zeta} $, more precisely,
$$
V_{x_0,\zeta}=\{x\in\zeta\colon \eta_{x,x_0}\in E^{\ssup{+}}_{x_0,\zeta}\}\qquad \mbox{ and } \;E_{x_0,\zeta}^{\ssup{\rm nbd}}=\{\eta_{x,y}\in E^{\ssup{\rm ext}}_{x_0,\zeta}\colon x,y\in V_{x_0,\zeta}\}.
$$

\begin{figure}[t]
\centering
\subfigure[$ E^{\ssup{\rm ext}}_{x_0,\zeta} $]{
\includegraphics[scale=1.25]{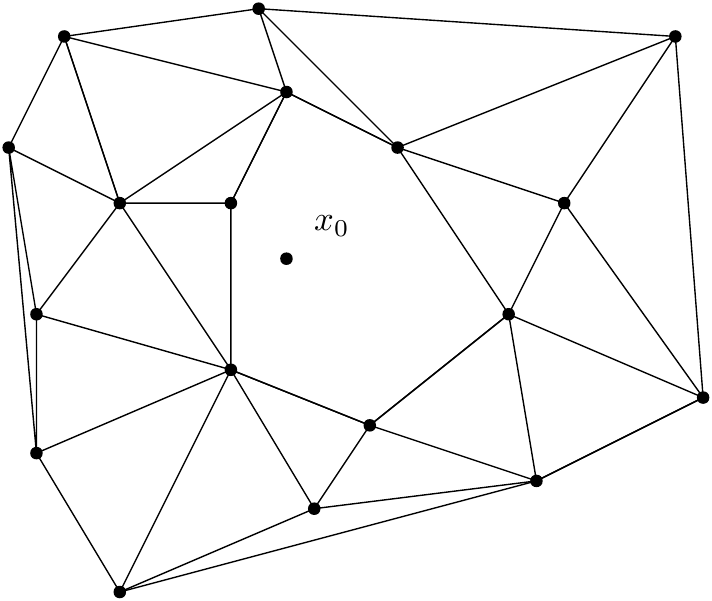}}
\subfigure[$ G_{x_0,\zeta} $]{
\includegraphics[scale=1.25]{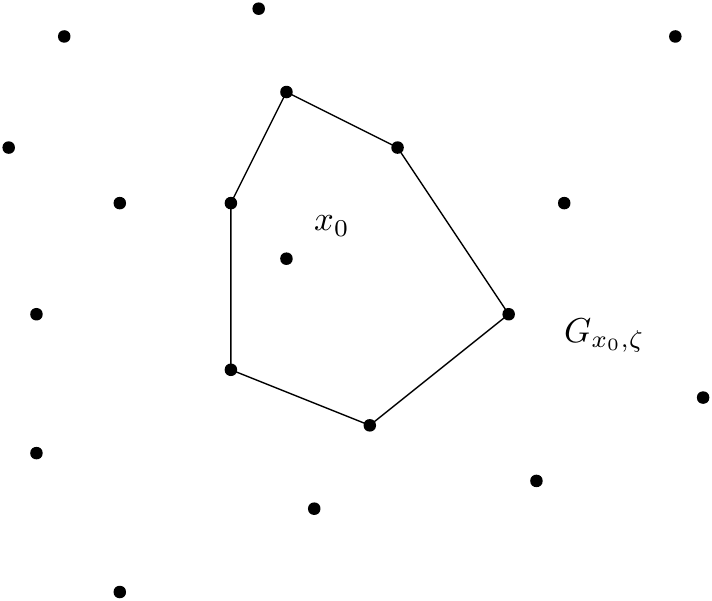}}
\caption{The exterior graph $ E^{\ssup{\rm ext}}_{x_0,\zeta} $ and the boundary graph $ G_{x_0,\zeta}$ }
\label{figboundarygraph}
\end{figure}

\bigskip

The graph $G_{x_0,\zeta} $ splits the plane into two regions. The region containing $x_0 $ is called the neighbourhood of $ x_0 $ whereas $ G_{x_0,\zeta} $ is called the boundary (graph) of the neighbourhood of $ x_0 $, see Figure~\ref{figboundarygraph}. Having the edge drawing mechanism $ \mu^{\ssup{\rm ext}}_{x_0,\zeta} $ define
\begin{equation}\label{edgem}
\mu_{\ssup{\rm ext},\zeta}^{\ssup{q}}(\d E):=\frac{q^{K(\zeta,E)}\mu_{x_0,\zeta}^{\ssup{\rm ext}}(\d E)}{\int\;q^{K(\zeta,E)}\mu_{x_0,\zeta}^{\ssup{\rm ext}}(\d E)}.
\end{equation}

The main task is the find an uniform upper bound (independent of $ \zeta $) for the expected number of connected components of $ (\zeta,E) $ that intersect the boundary $G_{x_0,\zeta} $, where $ E$ is sampled from $ \mu_{\ssup{\rm ext},\zeta}^{\ssup{q}} $. The bound will enable us to estimate  certain conditional  Papangelou intensities from below and above. This in turn  allows to dispense  background potentials used  in recent work (\cite{BBD04,AE16}).  Given $ \zeta\in\Omega_{\L,\omega} $ and $ x_0\in\L,x_0\notin\zeta $,  and the Delaunay graph $ \Del_2(\zeta) $ we denote $ N^{\ssup{\rm cc}}_{x_0}(\zeta,E) $ for any $ E\subset E_{x_0,\zeta}^{\ssup{\rm ext}} $ the number of connected components that intersect $ G_{x_0,\zeta} $.

\begin{theorem}[\textbf{Number of connected components}]\label{THMc}
Let $ \beta> q$. For all $ \L\Subset\R^2 $  there exists $ 0<\alpha=\alpha(R,q,\beta)= 1+6R^2\pi^2r^{-2}\big(1+\frac{2q\pi^2r^2}{3\beta}\big)  <\infty $ with $ r=1\wedge\frac{R\pi}{2} $ such that
$$
\int\;N^{\ssup{\rm cc}}_{x_0}(\zeta,E)\,\mu^{\ssup{q}}_{\ssup{\rm ext},\zeta}(\d E)\le \alpha
$$
for all  $ \omega\in\Omega $ and $ \zeta\in\Omega_{\L,\omega} $ and for all $ x_0\in\L $ with $ x_0\notin\zeta $. Note that, as $ \beta>q $, 
\begin{equation}\label{alphastar}
\alpha(R,q,\beta)\le \alpha^*(R,q):=1+6R^2\pi^2r^{-2}\big(1+\frac{2\pi^2r^2}{3}\big).
\end{equation}
Moreover,
\begin{equation}
\lim_{\beta\to\infty}\alpha(R,q,\beta)=1+6R^2\pi^2r^{-2}=\begin{cases} 25 & \mbox{ if } R<2/ \pi,\\
1+6R^2\pi^2 & \mbox{ if } R\ge 2/\pi  \,.
\end{cases}
\end{equation}
\end{theorem}

\begin{proofsect}{Proof} The proof is in Section~\ref{number}.\qed
\end{proofsect}

We need to study the change of $ K(\zeta,E), E\subset E_\zeta\subset\Escr $ when adding the point $ x_0\notin\zeta $. Adding a point $x_0 $ to $ \zeta $ without considering the change to $ E $ will always increase the number of connected components by one. On the other hand, the augmentation of a single edge $ \eta $ to $E$ can result in the connection of two different connected components, leaving one. Therefore,
\begin{equation}\label{change}
\begin{aligned}
K(\zeta\cup\{x_0\},E)-K(\zeta,E)&=1,\\
-1\le K(\zeta,E\cup\eta)-K(\zeta,E)&\le 0.
\end{aligned}
\end{equation}
 
\subsection{Delaunay Edge-Percolation}\label{tileperc}
 We establish the existence of edge percolation for the Delaunay random-cluster measure $ C_{\L,\om} $ when $ z$ and the parameter $ \beta $ are sufficiently large. Note that for any $ \Delta\Subset\R^2 $ we write
  $$
  N_{\Delta\leftrightarrow\infty}(\zeta,E)=\#\{x\in\zeta\cap\Delta\colon x \mbox{ belongs to an }  \infty-\mbox{cluster of } (\zeta,E\cap E_\zeta)\}.
  $$  
  The key step in our results is the following percolation result.
 \begin{proposition}\label{Prop-per}
  Suppose all the assumptions hold and that $ z $ and $ \beta $ are sufficiently large. Suppose that $ \L $ is a finite union of cells $ \Delta_{k,l} $ defined in \eqref{cell} Appendix A. Then there exists $\eps>0 $ such that
  $$
  \int\, N_{\Delta\leftrightarrow\infty}\,\d C_{\L,\om}\ge \eps
  $$ for any cell $\Delta=\Delta_{k,l} $, any finite union $ \L $ of cells and any admissible pseudo-periodic boundary condition $ \om \in\O_\L^*$.  
  \end{proposition}

\noindent\textbf{Proof of Proposition~\ref{Prop-per}}: We split the proof in several steps and Lemmata below.  
Our strategy to establish percolation in the Delaunay  random-cluster model is to compare it to mixed site-bond percolation on $ \Z^2 $ (see appendix~\ref{AppD} on mixed site-bond percolation). First we employ a coarse-graining strategy to relate each site $(k,l)\in\Z^2 $ to a cell which is a union of parallelotopes. The second step is to consider the links (bonds) of two good cells. In order to establish mixed site-bond percolation we need to define when cells are good (open) and when two neighbouring cells are linked once they are open which happens when the corresponding link (bond) is open as well. This link establishes an open connection in our Delaunay graph $ \Del_2(\zeta) $. We extend the coarse graining  method recently used in \cite{AE16}. 

\noindent \textbf{Step 1: Coarse graining.}

\begin{figure}[h!]
\caption{The $\ell$-partitioning of $ \L $.The shaded boxes are the link boxes}\label{figa}
\includegraphics[scale=1.22]{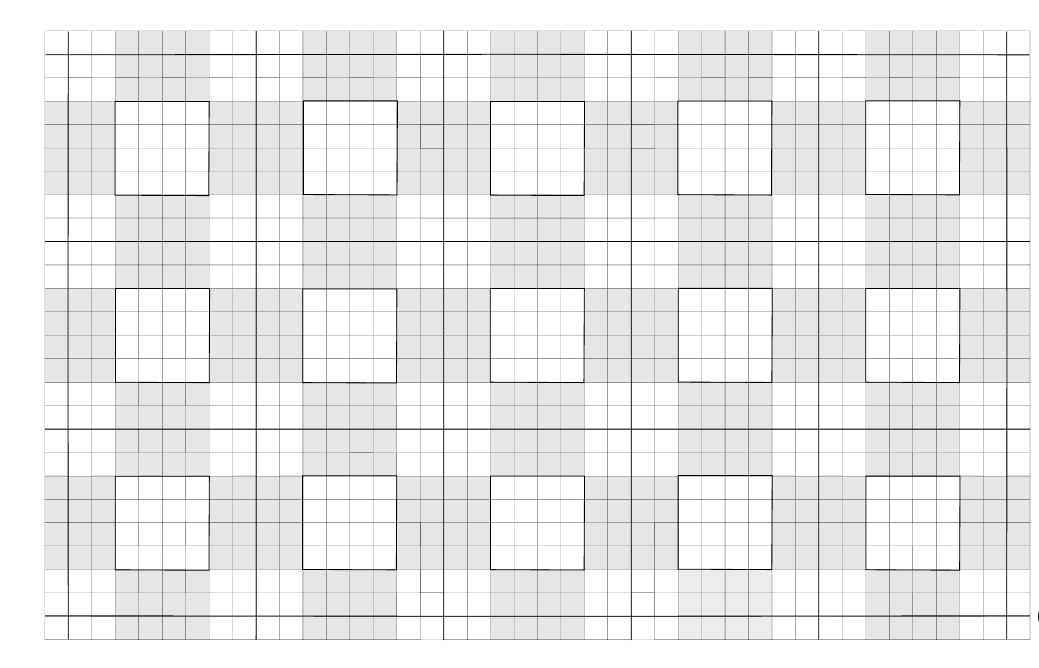}
\end{figure}

\begin{figure}[h!]
\caption{The shaded area is the union of the Voronoi cells with centre $\in H_{k:k+1,l}(\zeta) $}\label{figb}
\includegraphics[scale=1.52]{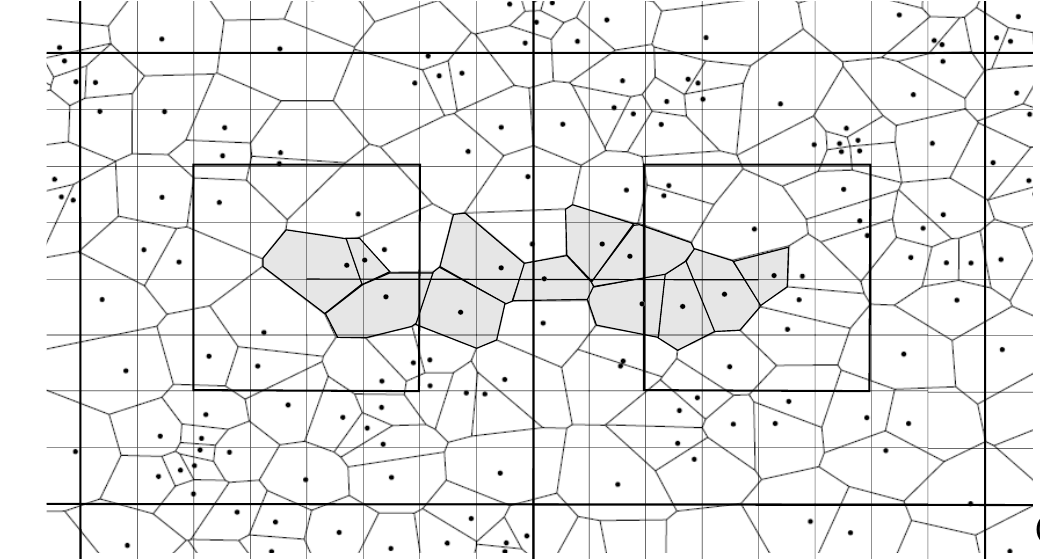}
\end{figure}

\noindent  Let $ \L=\L_n\subset\R^2 $ be the parallelotope  given as the  finite union of cells \eqref{cell} with side length $\ell $, i.e.,
 $$
 \L_n=\bigcup_{(k,l)\in \{-n,\ldots,n\}^2}\Delta_{k,l}\quad \mbox{ and } \Delta_{k,l}=\bigcup_{i,j=0}^7\Delta_{k,l}^{i,j},
 $$ 
 where $ \Delta_{k,l}^{i,j} $ are parallelotopes with side length $ \ell/8 $ and each parallelotope $ \Delta_{k,l} $ has side length $ \ell $ where the coordinate systems is the canonical one, that is, $ \Delta_{k,l}^{0,0} $ is the parallelotope in the bottom right corner.  The union of the $16 $ smaller parallelotopes towards the centre of $ \Delta_{k,l} $ is denoted
 $$
 \Delta_{k,l}^-=\bigcup_{i,j=2}^5\Delta_{k,l}^{i,j},
 $$
 see Figure~\ref{figa}. 
 These sites will act as the sites in the mixed site-bond percolation model on $ \Z^2 $. Finally, we define the link boxes between $ \Delta_{k,l}^- $ and $ \Delta_{k+1,l}^- $ as
 \begin{equation}\label{link}
 \Delta_{\ssup{\rm link}}^{k:k+1,l}=\Big\{\bigcup_{j=0}^3\Delta_{k,l}^{6,j+2}\Big\}\cup\Big\{\bigcup_{j=0}^3\Delta_{k,l}^{7,j+2}\Big\}\cup\Big\{\bigcup_{j=0}^3\Delta_{k+1,l}^{0,j+2}\Big\}\cup\Big\{\bigcup_{j=0}^3\Delta_{k+1,l}^{1,j+2}\Big\}
 \end{equation}
 which act as the bonds in the mixed site-bond percolation model on $\Z^2 $, see Figure~\ref{figa}. 
 We shall choose
 \begin{equation}\label{choiceoflength}
 \ell\in (0,\frac{R}{2\sqrt{3}}]
 \end{equation}
 to ensure that we can open edges in neighbouring parallelotopes. This completes the coarse graining  set-up. We establish percolation in the mixed site-bond percolation model on $ \Z^2 $, that is, the existence of an infinite chain of open sites and open bonds, and we relate it to the existence of an infinite connected component of open edges in $ \Del_2 $. This infinite connected component will connect with the complements of any finite boxes, and thus this connected component corresponds to an infinite connected component of edges where all sites carry the mark agreed for the boundary condition. To do this, we define  $ CB_{k:k+1,l} $ to be the straight line segment between the centres of the parallelotopes $ \Delta_{k,l} $ and $ \Delta_{k+1,l} $  and let
 \begin{equation}\label{horizontal}
 H_{k:k+1,l}(\zeta)=\{x\in\zeta\colon \Vor_\zeta(x)\cap CB_{k:k+1,l}\not=\emptyset\}
 \end{equation}
 be the subset of points of the configuration $ \zeta $, whose Voronoi cells intersect the line segment $ CB_{k:k+1,l} $, see Figures~\ref{figa} and \ref{figb}.
 
 \bigskip

  \bigskip

 \bigskip
 
We  need to consider the  distribution of the points  given by the marginal distribution
 $$
 M_{\L,\om}=C_{\L,\om}(\cdot\times\Escr)
 $$ of the Delaunay random-cluster measure on $ \O$. Note that  \eqref{DRCmeasure} can be written as
$$
 C_{\L,\om}(\d\zeta,\d E)=M_{\L,\om}(\d\zeta)\mu^{\ssup{q}}_{\zeta,\L}(\d E),\qquad \mu^{\ssup{q}}_{\zeta,\L}(\d E)=\frac{q^{K(\zeta,E)}\mu_{\zeta,\L}(\d E)}{\int\,q^{K(\zeta,E)}\mu_{\zeta,\L}(\d E)}.
$$
  We define $ h_\L $ to be the Radon-Nikodym density of $ M_{\L,\om} $ with respect to $ \Pi_{\L,\om}^z $, i.e., for $ \zeta\in\O_{\L,\om} $,
  $$
  h_\L(\zeta):=\Zscr_{\L}(\om)^{-1}\int\,q^{K(\zeta,E)}\mu_{\zeta,\L}(\d E).
  $$

\medskip 
 
 In the following lemma we derive a bound for the Papangelou conditional intensity of $ M_{\L,\zeta} $.

\begin{lemma}\label{L:Papangelou}
For any $ \L\Subset\R^2 $ and any admissible boundary condition $ \om\in\O_\L^* $ and $ M_{\L,\zeta}$-almost all $ \zeta\in\O_{\L,\om} $ and a point $ x_0\in\L $ with $ x_0\notin\zeta $,
\begin{equation}
\frac{h_\L(\zeta\cup\{x_0\})}{h_\L(\zeta)}\ge q^{1-\alpha},
\end{equation}
where $ \alpha\in (0,\infty) $ is given in Theorem~\ref{THMc}.
\end{lemma}

\begin{proofsect}{Proof of Lemma~\ref{L:Papangelou}}
Recall the different edge drawing mechanisms    $ \mu^{\ssup{-}}_{x_0,\zeta},\mu^{\ssup{+}}_{x_0,\zeta} $, and $ \mu^{\ssup{\rm ext}}_{x_0,\zeta} $ on $ E_{x_0,\zeta}^{\ssup{\rm ext}}, E^{\ssup{+}}_{x_0,\zeta} $, and $ E^{\ssup{-}}_{x_0,\zeta}$,  respectively, and the definition of the probability measure $ \mu_{\ssup{\rm ext},\zeta}^{\ssup{q}} $ in \eqref{edgem}. It follows that
$$
\begin{aligned}
\frac{h_\L(\zeta\cup\{x_0\})}{h_\L(\zeta)}&=\frac{\int\, q^{K(\zeta,E_2)}\int\,q^{K(\zeta\cup\{x_0\}, E_1\cup E_2)-K(\zeta,E_2)}\mu^{\ssup{+}}_{x_0,\zeta}(\d E_1)\mu^{\ssup{\rm ext}}_{x_0,\zeta}(\d E_2)}{\int\, q^{K(\zeta,E_3)}\int\,q^{K(\zeta,E_3\cup E_4)-K(\zeta,E_3)}\mu^{\ssup{-}}_{x_0,\zeta}(\d E_4)\mu^{\ssup{\rm ext}}_{x_0,\zeta}(\d E_3)}\\
&\ge \iint\, q^{K(\zeta\cup\{x_0\}, E_1\cup E_2)-K(\zeta,E_2)}\mu^{\ssup{+}}_{x_0,\zeta}(\d E_1)\mu^{\ssup{q}}_{\ssup{\rm ext},\zeta}(\d E_2)\\
&\ge q^{\iint (K(\zeta\cup\{x_0\}, E_1\cup E_2)-K(\zeta,E_2)) \mu^{\ssup{+}}_{x_0,\zeta}(\d E_1)\mu^{\ssup{q}}_{\ssup{\rm ext},\zeta}(\d E_2)},
\end{aligned}
$$
since
$$
K(\zeta, E_3\cup E_4)-K(\zeta,E_3)\le 0 $$ 
due to the fact that adding any edge from $E_4 \subset E_{x_0,\zeta}^{\ssup{-}}$ will only fuse connected components of the remaining graph. The second inequality is just Jensen's inequality applied to the convex function $ x\mapsto q^x $. Note that new edges from $ E_{x_0,\zeta}^{\ssup{+}} $, made by the insertion of $ x_0 $ to the configuration $ \zeta $, are edges connecting $ x_0 $ to points in $ \zeta $ and are open with respect to $ \mu^{\ssup{+}}_{x_0,\zeta} $, and therefore
\begin{equation}
K(\zeta\cup\{x_0\}, E_1\cup E_2)-K(\zeta,E_2)\ge - N^{\ssup{\rm cc}}_{x_0}(\zeta,E_2) +1
\end{equation} for any $ E_1\subset E_{x_0,\zeta}^{\ssup{+}} $ and any $ E_2\subset E_{x_0,\zeta}^{\ssup{\rm ext}} $,  and we conclude with the statement using Theorem~\ref{THMc}. 
 \qed
\end{proofsect}
An important component of our coarse graining method is to estimate the conditional probability that at least one point lies inside some $\Delta\subset\L $. For any $ \omega^\prime \in \Omega_{\Delta^{\rm c}} $  with $ \omega^\prime \cap \L^{\rm c}=\omega $ we denote by $ M_{\L,\Delta,\om^\prime} $ the conditional distribution of the configuration in $ \Delta $ given the configuration $ \omega^\prime $ in $ \Delta^{\rm c} $ relative to $M_{\L,\om} $.
The details of the construction of the regular conditional probability distribution can be found in \cite{E14} or \cite{GH96}. 
Having a uniform lower bound for the quotient $ h_\L(\zeta\cup\{x_0\})/h_\L(\zeta) $ allows to exhibit some control over the distribution  $ M_{\L,\nabla,\om^\prime} $ for any $ \nabla\subset\L $. In the following we write $ \nabla $ for any cell $ \Delta_{k,l}^{i,j}, i,j=0,\ldots,7 $. We fix
\begin{equation}
\eps=\frac{1-\sqrt{p_{\rm c}^{\ssup{\rm site}}(\Z^2)}}{4},
\end{equation}
where $ p_{\rm c}^{\ssup{\rm site}}(\Z^2) \in (0,1) $ is the critical probability for site percolation on $ \Z^2 $.

\begin{lemma}\label{lem1}
For $ \ell \in(0,\frac{R}{2\sqrt{3}}] $ here exists $ z_0=z_0(\alpha,q,\ell) $ such that  for all $ z> z_0 $ and    for all admissible pseudo-periodic boundary conditions $ \om \in\O_\L^*$,
$$
M_{\L,\nabla,\om^\prime}(\#\{\zeta\cap\nabla\}\ge 1)>1-\frac{\eps}{64}
$$
for all cells $ \nabla=\Delta_{k,l}^{i,j}, (k,l)\in\{-n,\ldots,n\}^2, i,j=0,\ldots, 7 $, and for any configuration $ \om^\prime\in\O_{\nabla^{\rm c}}  $ with $ \om^\prime\setminus\L=\om $.
\end{lemma}  
  
\begin{proofsect}{Proof}
Fix some $ \omega^\prime \in\Omega_{\nabla^{\rm c}} $ with $ \omega^\prime\setminus \L=\omega $. Then the statement follows immediately  from 
$$
\begin{aligned}
\frac{M_{\L,\nabla,\om^\prime}(\#\{\zeta\cap\nabla\}=1)}{M_{\L,\nabla,\om^\prime}(\#\{\zeta\cap\nabla\}=0)}=z\int_{\nabla}\frac{h_\L(\omega^\prime\cup\{x\})}{h_\L(\omega^\prime)}\,\d x\ge zq^{-\alpha}|\nabla|,
\end{aligned}
$$ where the inequality follows from Lemma~\ref{L:Papangelou}. It follows that
$$
M_{\L,\nabla,\om^\prime}(\#\{\zeta\cap\nabla\}=0)\le q^\alpha(z|\nabla|)^{-1},
$$ and hence 
$$
M_{\L,\nabla,\om^\prime}(\#\{\zeta\cap\nabla\}\ge 1)\ge 1-q^\alpha(z|\nabla|)^{-1}=1-q^\alpha(z\sqrt{3}/2\ell^2)^{-1}>1-\frac{\eps}{64}
$$ for $ z > z_0(\alpha,q,\ell)=\frac{2\cdot 64^2q^\alpha}{\eps\sqrt{3}\ell^2}>0$. Note that for all $ \ell\in(0,\frac{R}{2\sqrt{3}}] $,
\begin{equation}\label{boundz}
z_0(\alpha,q,\ell)\ge \frac{8\cdot 64^2\sqrt{3}q^\alpha}{\eps R^2}=:z_0^*(\alpha,q).
\end{equation}
\qed
\end{proofsect}

  \medskip

Let $ \Delta_{k,l} $ be an element of the  partitioning of $ \L $ and define $ F_{k,l}^{\ssup{\rm ext}} $ to be the event that each of the smaller boxes $ \Delta_{k,l}^{i,j} \subset\Delta_{k,l} $ that are not in the centre region $ \Delta_{k,l}^{\ssup{-}} $, contain at least one point. We call the elements in this event ``well-behaved'' configurations,
 \begin{equation}
 F_{k,l}^{\ssup{\rm ext}}=\bigcap_{\heap{i,j\in\{0,\ldots,7\}\colon}{\Delta_{k,l}^{i,j}\not\subset \Delta_{k,l}^{\ssup{-}}}}\big\{\zeta\in\O_{\L,\omega}\colon \#\{\zeta\cap\Delta_{k,l}^{i,j}\}\ge 1\big\}.
 \end{equation} 
  
  \begin{lemma}\label{L:PapangelouB}
For any $ \L\Subset\R^2 $ and any admissible boundary $ \om\in\O_\L^* $ and $ \ell$-partitioning of $ \L $  and well-behaved configuration $ \zeta\in F_{k,l}^{\ssup{\rm ext}} $ and a point $ x_0\in\Delta^{\ssup{-}}_{k,l} $ with $ x_0\notin\zeta $ for any $k,l\in\{-n,\ldots,n\} $,
\begin{equation}
\frac{h_\L(\zeta\cup\{x_0\})}{h_\L(\zeta)}\le q^{\alpha},
\end{equation}
where $ \alpha\in (0,\infty) $ is given in Theorem~\ref{THMc}.
\end{lemma}

\begin{proofsect}{Proof of Lemma~\ref{L:PapangelouB}}  
 Then, adapting similar steps in the proof of Lemma~\ref{L:Papangelou}, we see that
$$
 \begin{aligned}
\frac{h_\L(\zeta\cup\{x_0\})}{h_\L(\zeta)}&=\frac{\int\, q^{K(\zeta,E_2)}\int\,q^{K(\zeta\cup\{x_0\}, E_1\cup E_2)-K(\zeta,E_2)}\mu^{\ssup{+}}_{x_0,\zeta}(\d E_1)\mu^{\ssup{\rm ext}}_{x_0,\zeta}(\d E_2)}{\int\, q^{K(\zeta,E_3)}\int\,q^{K(\zeta,E_3\cup E_4)-K(\zeta,E_3)}\mu^{\ssup{-}}_{x_0,\zeta}(\d E_4)\mu^{\ssup{\rm ext}}_{x_0,\zeta}(\d E_3)}\\
&\le \frac{q \int\,q^{K(\zeta,E_2)}\,\mu^{\ssup{\rm ext}}_{x_0,\zeta}(\d E_2)}{ \int\, q^{K(\zeta,E_3)}\int\,q^{K(\zeta,E_3\cup E_4)-K(\zeta,E_3)}\mu^{\ssup{-}}_{x_0,\zeta}(\d E_3)\mu^{\ssup{\rm ext}}_{x_0,\zeta}(\d E_4)}   \\
&=  q\Big(\iint\,q^{K(\zeta,E_3\cup E_4)-K(\zeta,E_4)}\,\mu^{\ssup{-}}_{x_0,\zeta}(\d E_3)\mu^{\ssup{q}}_{\ssup{\rm ext},\zeta}(\d E_4)\Big)^{-1},
\end{aligned}
$$
where we used the inequality
\begin{equation}
K(\zeta\cup\{x_0\},E_1\cup E_2)-K(\zeta,E_2)\le 1, \qquad E_1\subset E^{\ssup{+}}_{x_0,\zeta}, E_2\subset E^{\ssup{\rm ext}}_{x_0,\zeta},
\end{equation}
as all new edges are connected to $ x_0 $ and thus can at most built one additional open component. We apply Jensen's inequality to the integral in the denominator to obtain the upper bound
\begin{equation}
\frac{h_\L(\zeta\cup\{x_0\})}{h_\L(\zeta)}\le q\Big(q^{\iint  (K(\zeta,E_3\cup E_4)-K(\zeta,E_4))\,\mu^{\ssup{-}}_{x_0,\zeta}(\d E_3)\mu^{\ssup{q}}_{\ssup{\rm ext},\zeta}(\d E_ 4)  }\Big)^{-1}.
\end{equation}
For all configurations $ \zeta\in F_{k,l}^{\ssup{\rm ext}} $ and $ x_0\in\Delta^{\ssup{-}}_{k,l} $ with $ x_0\not\in\zeta $ and all $ x\in V_{x_0,\zeta} $ we have that $ |x-x_0|<3\sqrt{3} \ell $. This ensures that
$$
V_{x_0,\zeta}\subset B_{3\sqrt{3}\ell}(x_0),
$$
where $ B_{3\sqrt{3}}\ell(x_0)  $ is the ball of radius $ 3\sqrt{3}\ell $ around $x_0 $. Therefore, since $x,y\in V_{x_0,\zeta} $ for all $ \eta_{x,y}\in E^{\ssup{-}}_{x_0,\zeta} $, it follows that adding edges in $ E^{\ssup{-}}_{x_0,\zeta} $ can only fuse together two connected components (reducing the number of connected components by one) if they each intersect $ V_{x_0,\zeta} $. Hence,
$$
K(\zeta,E_3\cup E_4)-K(\zeta,E_4)\ge -N^{\ssup{\rm cc}}_{x_0}(\zeta,E_4)+1, \quad E_3\subset E^{\ssup{-}}_{x_0,\zeta},
$$
and thus with Theorem~\ref{THMc} we conclude with the statement. 

  \qed
  \end{proofsect}
  
 We get a lower bound for the conditional probability, given well-behaved configurations $ \zeta\in F^{\ssup{\rm ext}}_{k,l} $, that $ \Delta^{\ssup{-}}_{k,l} $ contains no more that $ m=m(z) $ points of the well-behaved configurations $ \zeta $.

\begin{lemma}\label{L:mbound}
Given any $ \ell$-partitioning of $ \L $ with admissible boundary condition $ \omega\in\Omega $ and admissible boundary condition $\zeta^\prime\in\Omega $ for $ \Delta^{\ssup{-}}_{k,l} $ with $ \zeta^\prime\cap\L^{\rm c}=\omega $ and $ \zeta^\prime\in F^{\ssup{\rm ext}}_{k,l} $ for any $ k,l\in\{-n,\ldots,n\} $, the estimate
$$
M_{\L,\Delta^{\ssup{-}}_{k,l},\zeta^\prime}\big(\#\{\zeta\cap\Delta^{\ssup{-}}_{k,l}\}\le \floor{m(z)}\big)>1-\eps
$$ holds for $ m(z)=2\eps^{-1}q^{\alpha}|\Delta^{\ssup{-}}_{k,l}|z$.
\end{lemma}

\begin{proofsect}{Proof}
Note that we can write $M_{\L,\Delta^{\ssup{-}}_{k,l},\zeta^\prime}(\d\zeta) $ as
$$
M_{\L,\Delta^{\ssup{-}}_{k,l},\zeta^\prime}(\d\zeta)=\frac{1}{Z_{\L,\Delta^{\ssup{-}}_{k,l}}(\zeta^\prime)}h_\L(\zeta\cup\zeta^\prime)\,\Pi^z_{\Delta^{\ssup{-}}_{k,l}}(\d\zeta),
 $$ where $ Z_{\L,\Delta^{\ssup{-}}_{k,l}}(\zeta^\prime) $ is the normalisation. Using the well-known fact
 $$
 \int\,f(\zeta)\,\Pi^z_{\Delta^{\ssup{-}}_{k,l}}(\d\zeta)=\ex^{-z|\Delta^{\ssup{-}}_{k,l}|}\sum_{n=0}^\infty\frac{z^n}{n!}\int_{\Delta^{\ssup{-}}_{k,l}}\,f(\{x_1,\ldots,x_n\})\,\d x_1\cdots\d x_n
 $$ for any observable $ f $ of the underlying Poisson process and writing $ N(\zeta)=\#\{\zeta\cap\Delta^{\ssup{-}}_{k,l}\} $ we obtain, setting $Z^\prime=\ex^{-z|\Delta^{\ssup{-}}_{k,l}|}/Z_{\L,\Delta^{\ssup{-}}_{k,l}}(\zeta^\prime) $ for brevity, the probability for any $n\in\N $,
 $$
 \begin{aligned}
 &M_{\L,\Delta^{\ssup{-}}_{k,l},\zeta^\prime}(N =n+1)=\frac{z^{n+1}}{(n+1)!}Z^\prime\int_{\big(\Delta^{\ssup{-}}_{k,l}\big)^{n+1}} h_\L(\{x_1,\ldots,x_{n+1}\}\cup\zeta^\prime)\,\d x_1\cdots \d x_{n+1}\\
 &=\frac{z^n}{n!}\big(\frac{z}{n+1}\big)Z^\prime\int_{\big(\Delta^{\ssup{-}}_{k,l}\big)^n}\,\int_{\Delta^{\ssup{-}}_{k,l}}\, h_\L(\{x_1,\ldots,x_n\}\cup\zeta^\prime\cup\{x\})\,\d x\d x_1\cdots\d x_n\\
 &=\frac{z^n}{n!}\big(\frac{z}{n+1}\big)Z^\prime\int_{\big(\Delta^{\ssup{-}}_{k,l}\big)^n}\,h_\L(\{y\}\cup\zeta^\prime) g_{\Delta^{\ssup{-}}_{k,l},\zeta^\prime}(y)\, \d x\d y\end{aligned}
 $$
where 
$$
g_{\Delta^{\ssup{-}}_{k,l},\zeta^\prime}(y)=\int_{\Delta^{\ssup{-}}_{k,l}}\frac{h_\L(\{y\}\cup\zeta^\prime\cup\{x\})}{h_\L(\{y\}\cup\zeta^\prime)}\,\d x,  \quad y=\{x_1,\ldots,x_n\}.
$$
We obtain that
$$
\begin{aligned}
& M_{\L,\Delta^{\ssup{-}}_{k,l},\zeta^\prime}(N =n+1)=\big(\frac{z}{n+1}\big)\int_{\Omega_{\Delta^{\ssup{-}}_{k,l},\zeta^\prime}}\,g_{\Delta^{\ssup{-}}_{k,l},\zeta^\prime}(\zeta)\1\{N(\zeta)=n\}\,M_{\L,\Delta^{\ssup{-}}_{k,l},\zeta^\prime}(\d\zeta)=\big(\frac{z}{n+1}\big)\times\\
&\times M_{\L,\Delta^{\ssup{-}}_{k,l},\zeta^\prime}(N=n) \int_{\Omega_{\Delta^{\ssup{-}}_{k,l},\zeta^\prime}}\,g_{\Delta^{\ssup{-}}_{k,l},\zeta^\prime}(\zeta) M_{\L,\Delta^{\ssup{-}}_{k,l},\zeta^\prime}(\d\zeta|N=n).\end{aligned}
$$

By Lemma~\ref{L:PapangelouB} we can bound the function $ g_{\Delta^{\ssup{-}}_{k,l},\zeta^\prime}(\zeta) $ for point configurations $\zeta\cup\zeta^\prime\in F^{\ssup{\rm ext}}_{k,l} $ from above by $ q^\alpha|\Delta^{\ssup{-}}_{k,l}| $. Therefore,
$$
\begin{aligned}
\frac{M_{\L,\Delta^{\ssup{-}}_{k,l},\zeta^\prime}(N =n+1)}{M_{\L,\Delta^{\ssup{-}}_{k,l},\zeta^\prime}(N =n)}&=\Big(\frac{z}{n+1}\Big)\int_{\Omega_{\Delta^{\ssup{-}}_{k,l},\zeta^\prime}}\,  g_{\Delta^{\ssup{-}}_{k,l},\zeta^\prime}(\zeta) M_{\L,\Delta^{\ssup{-}}_{k,l},\zeta^\prime}(\d\zeta |N =n)\\
& \le\Big(\frac{z}{n+1}\Big)q^\alpha|\Delta^{\ssup{-}}_{k,l}|\int_{\Omega_{\Delta^{\ssup{-}}_{k,l},\zeta^\prime}}  M_{\L,\Delta^{\ssup{-}}_{k,l},\zeta^\prime}(\d\zeta|N =n)\le \frac{q^\alpha|\Delta^{\ssup{-}}_{k,l}| z}{n+1}.
\end{aligned}
$$
For all $n>m(z) $ we apply the previous step $ n-\floor{m(z)} $ times to obtain
$$
\begin{aligned}
 M_{\L,\Delta^{\ssup{-}}_{k,l},\zeta^\prime}(N =n)&\le \frac{q^\alpha|\Delta^{\ssup{-}}_{k,l}| z}{n}\,  M_{\L,\Delta^{\ssup{-}}_{k,l},\zeta^\prime}(\d\zeta|N =n-1)\\
 &\le  \frac{1}{n!}\big(q^\alpha|\Delta^{\ssup{-}}_{k,l}|z\big)^{n-\floor{m(z)}}\floor{m(z)}!.
 \end{aligned}
$$
It follows that
$$
\begin{aligned}
 M_{\L,\Delta^{\ssup{-}}_{k,l},\zeta^\prime}(N>m(z))&\le\sum_{n=\floor{m(z)}+1}^\infty\frac{1}{n!} \big(q^\alpha|\Delta^{\ssup{-}}_{k,l}|z\big)^{n-\floor{m(z)}}\floor{m(z)}!\\ 
 &\le \sum_{n=\floor{m(z)}+1}^\infty\Big(\frac{q^\alpha|\Delta^{\ssup{-}}_{k,l}|z}{\floor{m(z)}}\Big)^{n-\floor{m(z)}} ,
\end{aligned}
$$ and thus
$$
 M_{\L,\Delta^{\ssup{-}}_{k,l},\zeta^\prime}(N>m(z))\le  \sum_{n=\floor{m(z)}+1}^\infty\big(\frac{\eps}{2}\big)^{n-\floor{m(z)}}=\sum_{n=1}^\infty\big(\frac{\eps}{2}\big)^n.
 $$
Since $ \eps<\frac{1}{2} $, the right hand side of the previous inequality is less than $ \eps $. We conclude with the statement.
\qed
\end{proofsect}

  \noindent \textbf{Step 2: Random-cluster measure $\widetilde{C}_{\L,\om}$.}
 We find a measure $ \widetilde{C}_{\L,\om} $ which is stochastically smaller than $ C_{\L,\om} $. Then using coarse graining  and comparison to mixed site-bond percolation on $ \Z^2 $ we establish percolation for $ \widetilde{C}_{\L,\om} $. Percolation for $ \widetilde{C}_{\L,\om} $ then implies percolation for the original random cluster measure $ C_{\L,\om} $. We base the definition of the measure $ \widetilde{C}_{\L,\om} $ on a coarse graining method originally introduced in \cite{Haggstrom00} and later extended and adapted  in \cite{AE16}.

 \noindent For given $ \zeta\in\O_{\L,\om} $ and $ \ell \in (0,\frac{R}{2\sqrt{3}}] $ let 
 $$
 \Del_2^*(\zeta)=\big\{\eta\in\Del_2(\zeta)\colon \phi_\beta(\ell(\eta))\ge g(\beta):=\log(1+\beta/64\ell^4)\big\},
 $$
 and let $ \widetilde{\mu}_{\zeta,\L} $ be the distribution of the random set $\{\eta\in E_\zeta\colon \upsilon(\eta)=1\} $  with
 \begin{equation}\label{siteperc}
 \P(\upsilon(\eta)=1)=\widetilde{p}_\ell(\eta)=\frac{1-\exp\{-g(\beta)\}}{1+(q-1)\exp\{-g(\beta)\}}\1_{\Del_2^*(\zeta)}(\eta).
 \end{equation}
 Note that $ \widetilde{\mu}_{\zeta,\L} $ depends on $ \L $ only via the configuration $ \zeta\in\O_{\L,\om}$. 
 Note also the important fact that $ \widetilde{p}_\ell $  is increasing in $\beta$, although to reduce excessive notation, we don't explicitly write this. 
 It is easy to show that $ \mu_{\zeta,\L}\succcurlyeq\widetilde{\mu}_{\zeta} $ by noting that
 $$
 \begin{aligned}
 \frac{p(\eta)}{q(1-p(\eta))}&=\frac{1-\exp(-\phi_\beta(\eta))}{q\exp(-\phi_\beta(\eta))}\ge\frac{1-\exp[-g(\beta)\1_{\Del^*_2(\zeta)}(\eta)]}{q\exp[-g(\beta)\1_{\Del^*_2(\zeta)}(\eta)]}= \frac{\widetilde{p}_\ell(\eta)}{q(1-\widetilde{p}_\ell(\eta))}, \qquad \eta\in\Del_2(\zeta),
 \end{aligned}
 $$ and using \cite{G94}.
 Hence, $ C_{\L,\om}\succcurlyeq\widetilde{C}_{\L,\om} $. As site percolation implies edge percolation we consider site percolation given by \eqref{siteperc}, that is, we open vertices in $ \Del_1(\zeta) $  independently of each other with probability $ \widetilde{p}_\ell$. Formally this is defined as follows.  We let $ \widetilde{\lambda}_{\zeta,\L} $ be the distribution of the random mark vector $ \widetilde{\sigma}\in E^\zeta $ where $ (\widetilde{\sigma}_x)_{x\in\zeta} $ are Bernoulli random variables satisfying
 \begin{equation}\label{ptilde}
 \begin{aligned}
 \P(\widetilde{\sigma}_x=1)&=\widetilde{p}_\ell\1_{\Del_1^*(\zeta)}(x)\\
 \P(\widetilde{\sigma}_x\not= 1)&=1-\widetilde{p}_\ell\1_{\Del_1^*(\zeta)}(x),
 \end{aligned}
 \end{equation}
 where $ \widetilde{p}_\ell $ is given in \eqref{siteperc} and $ \Del_1^*(\zeta) $ is the set of points that build the edges of $ \Del_2^*(\zeta) $.
Then the site percolation process is defined by the measure
\begin{equation}\label{sitecluster}
 \widetilde{C}^{\ssup{\rm site}}_{\L,\om}(\d\bzeta)=M_{\L,\om}(\d\zeta)\widetilde{\lambda}_{\zeta,\L}(\d\bzeta).
 \end{equation}\label{boundp}
 Note that for all $ \eta\in\Del_2^*(\zeta) $ and all $ \ell\in (0,\frac{R}{2\sqrt{3}}] $,
 \begin{equation}
 \widetilde{p}_\ell\ge \frac{1}{\frac{4qR^4}{q\beta} +1}=:\widetilde{p}.
 \end{equation}
 Note that $ \widetilde{p}=\widetilde{p}(\beta) $ is increasing in $ \beta $.

 \noindent\textbf{Step 3: Site-bond percolation.}
 
We now establish percolation for the random-cluster measure $ \widetilde{C}_{\L,\om} $.
 
\begin{lemma}\label{LemPercolation}
  There is a  $c>0 $ such that 
$$
\widetilde{C}^{\ssup{\rm site}}_{\L,\om}(\{\Delta\longleftrightarrow\L^{\rm c}\})\ge c>0
$$ for any $ \Delta=\Delta_{k,l}\subset\L, (k,l)\in\{-n,\ldots,n\}^2$, and any pseudo-periodic admissible boundary condition $ \om \in\O^*_{\L}$.
\end{lemma}

\begin{proofsect}{Proof of Lemma~\ref{LemPercolation}}

\noindent\textbf{Step (i) Probability that small cells have at least one point:} In the following we write $ \Delta=\Delta_{k,l} $  and $ \nabla=\Delta_{k,l}^{i,j} $ for any $k,l\in\{-n,\ldots,n\} $ and for any $ i,j\in\{0,1,\ldots,7\} $. For all configurations $ \zeta^\prime\in\Omega_{\Delta^{\rm c}} $ with $ \zeta^\prime\cap\L^{\rm c}=\omega $ we obtain with Lemma~\ref{lem1} the estimate
\begin{equation}\label{step1est}
M_{\L,\Delta,\zeta^\prime}(F^{\ssup{\rm ext}}_{k,l})\ge 1-\sum_{i,j\colon \Delta_{k,l}^{i,j}\not\subset\Delta^{\ssup{-}}_{k,l}}M_{\L,\Delta,\zeta^\prime}(\#\{\zeta\cap\Delta_{k,l}^{i,j}\}=0)<1-\frac{48\eps}{64}=1-\frac{3\eps}{4}.
\end{equation}
Define the following two events, first the event
$$ 
G_{k,l}=\{\zeta\in\O_{\L,\om}\colon\#\{\zeta\cap\Delta_{k,l}^{\ssup{-}}\}\le m(z)\}
$$ that there are at most $ m(z) $ points in $ \zeta $ in the centre $ \Delta_{k,l}^{\ssup{-}} $ of $ \Delta_{k,l} $, and the event that all smaller cells in that centre contain at least one point,
$$
F^{\ssup{-}}_{k,l}=\bigcap_{i,j\colon \Delta_{k,l}^{i,j}\subset\Delta_{k,l}^{\ssup{-}}}\big\{\zeta\in\O_{\L,\om}\colon\#\{\zeta\cap \Delta^{i,j}_{k,l}\}\ge 1\big\}.
$$
We have tacitly replaced $ m(z) $ by $ m(z)\vee 16 $. Both events depend on point configurations in the centre region $ \Delta^{\ssup{-}}_{k,l} $, and it suffices to estimate the probability of the intersection of these two events for any boundary condition outside of $ \Delta_{k,l} $ and any point configuration $ \zeta^{\dprime} $ in $ \Delta_{k,l}\setminus\Delta^{\ssup{-}}_{k,l} $. For all boundary conditions $ \zeta^{\sprime}\in\O_{\Delta_{k,l}^{\rm c}} $ with $ \zeta^{\sprime}\cap\L^{\rm c}=\om $,
$$
\begin{aligned}
&M_{\L,\Delta,\zeta^{\sprime}}(F^{\ssup{-}}_{k,l}\cap G_{k,l}\cap F^{\ssup{\rm ext}}_{k,l})=\\
&\int\,\1_{F_{k,l}^{\ssup{\rm ext}}}(\zeta^{\dprime}\cup\zeta^{\sprime})\Big[\int\,\1_{F_{k,l}^{\ssup{-}}}(\zeta)\1_{G_{k,l}}(\zeta)M_{\L,\Delta,\zeta^{\sprime}}(\d\zeta|\zeta=\zeta^{\dprime}\mbox{ on }\Delta\setminus\Delta^{\ssup{-}})\Big]M_{\L,\Delta,\zeta^{\sprime}}(\d\zeta^{\dprime}).
\end{aligned}
$$
Using Lemma~\ref{L:mbound} it follows that
\begin{equation}
\begin{aligned}
M_{\L,\Delta^{\ssup{-}}_{k,l},\zeta^{\dprime}\cup\zeta^{\sprime}}(F^{\ssup{-}}_{k,l}\cap G_{k,l})&\ge 1- M_{\L,\Delta^{\ssup{-}}_{k,l},\zeta^{\dprime}\cup\zeta^{\sprime}}((F^{\ssup{-}}_{k,l})^{\rm c})- M_{\L,\Delta^{\ssup{-}}_{k,l},\zeta^{\dprime}\cup\zeta^{\sprime}}((G_{k,l})^{\rm c})\\
&>1-\frac{16\eps}{64}-\eps=1-\frac{5\eps}{4},
\end{aligned}
\end{equation}
and hence we conclude with \eqref{step1est} that
\begin{equation}
\begin{aligned}
M_{\L,\Delta_{k,l},\zeta^{\sprime}}(F^{\ssup{-}}_{k,l}\cap G_{k,l}\cap F^{\ssup{\rm ext}}_{k,l})&>\big(1-\frac{5\eps}{4}\big)\int\,\1_{F_{k,l}^{\ssup{\rm ext}}}(\zeta^{\dprime}\cup\zeta^{\sprime}) M_{\L,\Delta,\zeta^{\sprime}}(\d\zeta^{\dprime})\\
&>\big(1-\frac{5\eps}{4}\big)\big(1-\frac{3\eps}{4}\big)>1-2\eps.
\end{aligned}
\end{equation}
\medskip

\noindent\textbf{Step (ii): Good cells and site percolation:}\\[1ex]
After these preparation steps we shall define when a cell $ \Delta_{k,l} $ is good. A good cell at $ (k,l) $ will result in the site $ (k,l)\in\Z^2 $ being open.  The next step is therefore to condition the marks of the points, that is, we pick $ (k,l)\in\{-n,\ldots,n\}^2 $ and consider the event $ C_{k,l} $ that each cell $ \Delta_{k,l}^{i,j} $ has at least one point, $\Delta^{\ssup{-}}_{k,l} $ contains no more than $ m(z) $ points   and all points in $ \Delta_{k,l}^{\ssup{-}}\cap\Del_1^*(\zeta) $ are carrying mark $1$,
$$
C_{k,l}=\{\bzeta\in\bO_{\L,\bo}\colon \zeta\in F_{k,l}^{\ssup{-}}\cap F^{\ssup{\rm ext}}_{k,l}\cap G_{k,l}  \mbox{ and } \sigma_{\bzeta}(x)=1 \mbox{ for all } x\in\Delta_{k,l}^{\ssup{-}}\cap\Del_1^*(\zeta)\}.
$$
A cell $ \Delta_{k,l} $ is declared to be ``good'' or ``open'' if $ C_{k,l} $ occurs.
Each vertex $ x\in\Del_1^*(\zeta) $ is open with probability $ \widetilde{p}_\ell $ (see \eqref{ptilde}). It follows that
$$
\begin{aligned}
\widetilde{C}_{\L,\om}^{\ssup{\rm site}}(C_{k,l})&\ge \int\,M_{\L,\Delta_{k,l},\zeta^{\sprime}}(\d\zeta)\1_{F^{\ssup{-}}_{k,l}}(\zeta)\1_{F^{\ssup{\rm ext}}_{k,l}}(\zeta)\1_{G_{k,l}}(\zeta)\widetilde{p}_\ell^{\#\{\Del_1^*(\zeta)\cap\Delta^{\ssup{-}}_{k,l}\}}\\&\ge \widetilde{p}_\ell^{\floor{m(z)}}M_{\L,\Delta,\zeta^{\sprime}}(F^{\ssup{-}}_{k,l}\cap G_{k,l}\cap F^{\ssup{\rm ext}}_{k,l}).
\end{aligned}
$$
Recall from Lemma~\ref{lem1} that there is $ z_0=z_0(\alpha,q,\ell) $ such that the estimates hold for all $ z\ge z_0 $.  Recall that $ |\Delta^{-}|=\frac{\sqrt{3}}{8}\ell^2\le\frac{R^2}{32\sqrt{3}} $ and $ \floor{m(z)}\le \frac{2q^\alpha R^2z}{32\sqrt{3}\eps} $. 
For all $ z\ge z_0(\alpha,q,\ell)\ge z_0^*(\alpha,q) $ choose $ \beta_0=\beta_0(q,R,z) >0  $ such that
\begin{equation}\label{beta0}
(\widetilde{p})^{\frac{q^{\alpha^* }R^2z}{16\sqrt{3}\eps}}\ge (1-2\eps)\quad \mbox{ for all } \beta\ge \beta_0\vee q,
\end{equation}
where $ \alpha^* $ is the bound for $ \alpha $ \eqref{alphastar} for any $ \beta>q $.
Then, for all $ \ell\in (0,\frac{R}{2\sqrt{3}}] $ and all $ z\ge z_0(\alpha^*,q,\ell)\ge z_0^*(\alpha^*,q) $ and all $ \beta\ge \beta_0\vee q $,
\begin{equation}
\widetilde{p}_\ell^{\floor{m(z)}}\ge \widetilde{p}^{\floor{m(z)}}\ge (1-2\eps).
\end{equation}
Combining the above estimates, we conclude, for all $ \beta\ge \beta_0 \vee q$, that
\begin{equation}\label{site}
\widetilde{C}^{\ssup{\rm site}}_{\Delta_{k,l},\zeta^{\sprime}}(C_{k,l})\ge (1-2\eps)^2>1-4\eps>(p_{\rm c}^{\ssup{\rm site}}(\Z^2))^{1/2}.
\end{equation}

\medskip

\noindent\textbf{Step (iii) Neighbouring good cells and link percolation:}\\[1ex]

If $\zeta\in C_{k,l} $, we say that the cell $ \Delta_{k,l} $ is a ``good'' cell. Two neighbouring cells $ \Delta_{k,l} $ and $ \Delta_{k+1,l} $ are said to be ``linked'' if the box $ \Delta_{\ssup{\rm link}}:=\Delta^{k:k+1,l}_{\ssup{\rm link}} $ defined in \eqref{link} has an intersection with $ \Del_1^*(\zeta) $ that contains only points of mark $1$. More precisely, the event that $ \Delta_{k,l} $ and $ \Delta_{k+1,l} $ are linked, is
$$
L^l_{k:k+1}=\{\bzeta\in\bO_{\L,\bo}\colon \sigma_{\bzeta}(x)=1 \mbox{ for all } x\in\Delta^{k:k+1,l}_{\ssup{\rm link}}\cap\Del_1^*(\zeta)\}.
$$
We also define
$$
F_{\ssup{\rm link}}=\Big(F_{k,l}^{\ssup{-}}\cap F_{k,l}^{\ssup{\rm ext}}\Big)\cap\Big(F^{\ssup{-}}_{k+1,l}\cap F^{\ssup{\rm ext}}_{k+1,l}\Big)
$$
and
$$
G_{\ssup{\rm link}}=\{\zeta\in\O\colon\#\{\zeta\cap\Delta_{\ssup{\rm link}}\}\le m(z)\},
$$
and let $\zeta^{\sprime} \in\O_{\Delta_{\ssup{\rm link}}^{\rm c}} $ be the boundary condition outside $ \Delta_{\ssup{\rm link}} $ such that $ \zeta^{\sprime}\cap\L^{\rm c}=\om $. The conditional probability that the cells $ \Delta_{k,l} $ and $ \Delta_{k+1,l} $ are linked, given they are both ``good'' cells, is then given by
\begin{equation}\label{linkp}
\begin{aligned}
\widetilde{C}^{\ssup{\rm site}}_{\Delta_{\ssup{\rm link}},\zeta^{\sprime}}(L^l_{k,k+1}|C_{k,l}\cap C_{k+1,l})&\ge \int\,\1_{G_{\ssup{\rm link}}}(\zeta)\widetilde{p}^{\floor{m(z)}}\,M_{\L,\Delta_{\ssup{\rm link}},\zeta^{\sprime}}(\d\zeta |F_{\ssup{\rm link}})\\
&\ge (1-\eps)(1-2\eps)\ge (1-4\eps)\ge  (p_{\rm c}^{\ssup{\rm site}}(\Z^2))^{1/2},
\end{aligned}
\end{equation}
where the second inequality comes from an adaptation of Lemma~\ref{L:mbound}, where $ \Delta_{\ssup{\rm link}} $ takes the role of $ \Delta^{\ssup{-}}_{k,l} $.  Then by \eqref{site}, \eqref{linkp} and the results of McDiarmid and Hammersley, in particular, \eqref{ineqmixed}, mixed site-bond percolation in $ \Z^2 $ occurs. There exists a chain of good cells  joined by open links from $ \Delta_{k,l}\subset\L $ to $ \L^{\rm c} $.  

It remains to check that the preceding  result implies $ \{\Delta\leftrightarrow \L^{\rm c}\} $ in the Delaunay graph. For this, we recall the set $ H_{k:k+1,l}(\zeta) $ in  \eqref{horizontal}. We know by construction that all edges $ \eta=\{x,y\}\in\Del_2(\zeta) $ that have a non-empty intersection with $ H_{k:k+1,l}(\zeta) $ satisfy $ |x-y|<2\sqrt{3}\frac{\ell}{8} $. This implies that $ H_{k:k+1,l}(\zeta)\subset \Del_1^*(\zeta) $. Let $x,y\in\zeta $ be such that $ \Vor_\zeta(x) $ and $ \Vor_\zeta(y) $ contain the centres of the boxes $ \Delta_{k,l} $ and $ \Delta_{k+1,l} $ respectively. Since $ H_{k:k+1,l}(\zeta)\subset\Del_1^*(\zeta) $, we can connect $x $ and $ y $ in the graph $ \Del_2^*(\zeta) $ inside $ \Delta^{\ssup{-}}_{k,l}\cup\Delta_{\ssup{\rm link}}\cup\Delta^{\ssup{-}}_{k+1,l} $. Hence, by  \eqref{site} and  \eqref{linkp}, there exits $c>0 $ such that  the following uniform lower bound holds
$$
\widetilde{C}^{\ssup{\rm site}}_{\L,\om}(\{\Delta\leftrightarrow\L^{\rm c}\})>c>0,
$$
and the proof of Lemma~\ref{LemPercolation} is finished.

\qed
\end{proofsect}

\medskip

\noindent\textbf{Step 4:  Finish of the proof of Proposition~\ref{Prop-per}}. The proof follows immediately from all previous steps as percolation in the site percolation measure $ \widetilde{C}^{\ssup{\rm site}}_{\L,\omega} $ implies percolation in the Delaunay random cluster measure $ C_{\L,\omega} $ due to stochastic dominance, 
$$C_{\L,\om}\succcurlyeq\widetilde{C}_{\L,\om}. $$

\qed

\subsection{Symmetry breaking of the mark distribution}\label{phase}

To relate the influence of the boundary condition on the mark of a single point  to the connectivity probabilities in the random-cluster model we follow \cite{GH96}. For any $ \Delta\subset\L $, $s\in E $, $\bzeta\in\bO_{\L,\bo} $ and $ (\zeta,E) $, with $ E\subset E_\zeta $ we define
 $$
 N_{\Delta,s}(\bzeta)=\#\{\zeta^{\ssup{s}}\cap\Delta\}.
 $$ 
 Then 
 $$
 \begin{aligned}
 N_{\Delta\leftrightarrow\L^{\rm c}}(\zeta,E)=\#\big\{ &x\in\zeta\cap\Delta\colon x\mbox{ belongs to a cluster connected to } \L^{\rm c}  \mbox{ in } E\cap\Del_2(\zeta)\big\}
 \end{aligned}
 $$ is the number of points in $ \zeta\cap\Delta $ connected to any point in $ \L^{\rm c} $ in the random graph $E\cap\Del_2(\zeta) $.
 Because of the edge-drawing mechanism, $ \{\Delta\leftrightarrow\L^{\rm c}\} =\{N_{\Delta\leftrightarrow\L^{\rm c}}>0\} $ is also the event that there exists a point in $ \zeta\cap\Delta $ connected to infinity in $ E\cap\Del_2(\zeta) $.
 
 The next Proposition is the key argument why percolation for the random cluster measures leads to a break of symmetry in the mark distribution.
 \begin{proposition}\label{Prop-sym}
 For any measurable $\Delta\subset\L \Subset\R^2$,
 $$
 \int\,(qN_{\Delta,1}-N_\Delta)\,\d\gamma_{\L,\bo} =(q-1)\int\,N_{\Delta\leftrightarrow\L^{\rm c}}\,\d C_{\L,\om}.
 $$
 \end{proposition}
 \begin{proofsect}{Proof}
 This is proved in \cite[Lemma~2.17]{E14} following ideas in  \cite{GH96}.
 \qed
 \end{proofsect}

\section{Number of connected components}\label{number}
In this section we are going to prove the main technical tool for our phase transition proof, namely, Theorem~\ref{THMc}. The proof of Theorem~\ref{THMc} is rather long, so we first outline the strategy. We want to bound the number of connected components in the graph $ E_{x_0,\zeta}^{\ssup{\rm ext}} $ that intersect the boundary graph $ G_{x_0,\zeta} $ under the edge drawing mechanism $ \mu^{\ssup{q}}_{\ssup{\rm ext},\zeta} $. We  also define $ G_B $ for the contraction of $ G_{x_0,\zeta} $ to a ball $B$ around $x_0 $, that is $ G_B=(V_B,E_B) $ with $ V_B=V_{x_0,\zeta}\cap B $ and $ E_B=\{\eta_{x,y}\in\Del_2(\zeta\cap B)\colon x,y\in V_B\} $. In the following we choose $ B=B_R(x_0) $ as due to the finite range condition the point $ x_0 $ cannot be connected to any point further away than $ R>0 $.
The pivotal point of the whole proof is to find an upper bound for the number of edges in the edge set $E_B $ that have length greater than some fixed real number. This allows us to construct families consisting of edges in $E_B $, defined by edge length, to balance the unbounded number of points against the increased probability that they are connected. The shorter the edge length, the greater the possible number of edges in the subset, but also the greater the probability that these edges are open. It turns out that such an upper bound can be found in the scenario where there are no defects in the geometry. These defects which we give the logical name ``kinks'' are defined below in Section~\ref{notation}. An upper bound cannot be found if the geometry of the contracted graph contains kinks, so we devise a plan to discount them. 

For $ R>0 $ the following three cases depend on the configuration $ \zeta $ and the point $ x_0 $. (i) $V_B  =\emptyset $ in which case there is no connection to any connected component of $ E_{x_0,\zeta}^{\ssup{\rm ext}} $, (ii) $V_B\supset V_{x_0,\zeta} $, and (iii) $V_B\cap V_{x_0,\zeta}\not=\emptyset $ and $ B_R(x_0)^{\rm c}\cap V_{x_0,\zeta}\not=\emptyset $. In case (ii) we have that $ E_B \subset E^{\ssup{\rm ext}}_{x_0,\zeta} $ but this does not hold in case (iii). This creates a problem when dealing with our edge drawing mechanism on $ E^{\ssup{\rm ext}}_{x_0,\zeta} $. To overcome this, we introduce an edge drawing mechanism on the contracted graph and build a structure that will allow us to compare events between the two different probability spaces. All our techniques rely heavily on geometric properties of the Delaunay tessellation.  

\subsection{Notation and geometric facts}\label{notation}
We introduce a polar coordinate system in $ \R^2 $ with $ x_0 $ being the pole, and we denote $L$ the polar axis in horizontal direction. For any $z\in\R^2 $, denote $\hat{z} $ be the angular coordinate of $z$ taken counter clockwise from the axis $L$. Given two points $x,y\in\R^2 $, $\overleftrightarrow{xy} $ denotes the unique straight line that intersects $x $ and $ y $  in the plane, $\overleftarrow{xy} $ denotes the half line that stops at $y$ and $  \overline{xy} $ denotes the line segment between $x $ and $ y $ only. Given two straight lines $ \ell_1,\ell_2 \subset\R^2 $ that intersect a point $z\in\R^2 $, $ \angle{(\ell_1,\ell_2)} $ denotes the angle between them. More precisely, it is the angle in order to rotate $ \ell_1 $ onto $ \ell_2 $ with $z$ being the centre of rotation. Notice that $ \angle{(\ell_1,\ell_2)}=\angle{(\ell_2,\ell_1)} $ only if $ \angle(\ell_1,\ell_2)=\pi/2$, however, it holds that $ \angle{(\ell_1,\ell_2)}+\angle{(\ell_2,\ell_1)}=\pi $. When we consider a triangle in the plane with vertices $ x,y,z $,  we refer to the interior angle at $y $ as $ \widehat{xyz} $. In this case, as we specify the interior angle, $ \widehat{xyz}=\widehat{zyx}$.

\begin{figure}[h!]
\caption{From top to bottom we have: (1.) A collection of points that neighbour $x_0 $ in the Delaunay (Voronoi) tessellation. (2.) A spoked chain $ \Gamma $ with pole $x_0$. (3.) The induced polygon $ P(\Gamma,x_0)$.}\label{figd}
\includegraphics[scale=1.0]{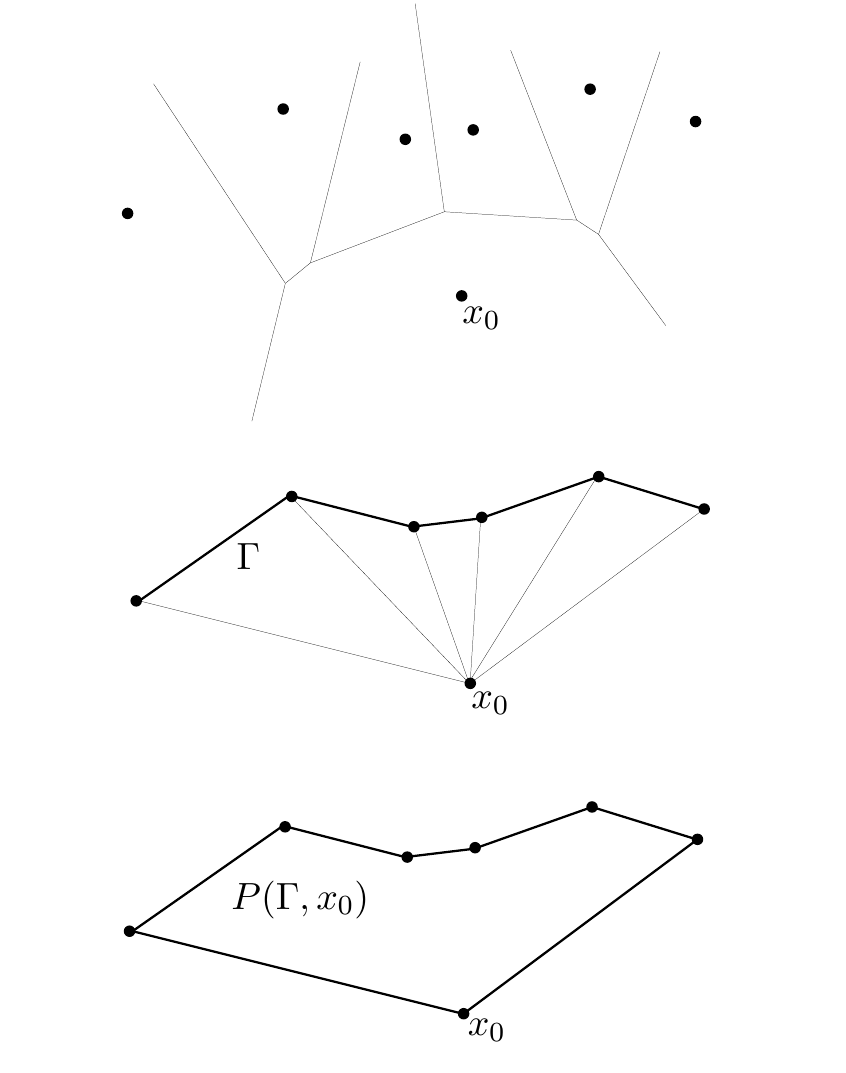}
\end{figure}

Given a set of points $V=\{x_i\in\R^2\colon 1\le i\le n\} $ with $ \hat{x}_1<\cdots<\hat{x}_n $, the graph
$$
\Gamma=\Big(V,\bigcup_{i=1}^{n-1}\eta_{x_i,x_{i+1}}\Big)
$$
is called a \textbf{spoked chain} with pole $ x_0 $ if $ \eta_{x_0,x_i}\in\Del_2(V\cup\{x_0\}) $ for all $ 1\le i\le n $. The polygon $ P(\Gamma,x_0) $ created by adding the point $ x_0 $ and edges $ \eta_{x_0,x_1} $ and $ \eta_{x_n,x_n} $ to $ \Gamma $ is called the induced polygon of $ \Gamma $, see Figure~\ref{figd}.

\bigskip

\bigskip

In order to quantify the number of connected components that intersect $V_B $, we analyse the shape of the contracted boundary graph  $ G_B $. For that we split $B$ into four quadrants, $Q_i\subset\R^2, i=1,2,3,4 $, where
$$
Q_i=\big\{z\in B\colon \frac{\pi}{2}(i-1)\le \hat{z}<\frac{\pi}{2}i\big\}.
$$
This allows to bound the number of connected components in one quadrant, and then  the final bound is just four times this bound. The reasons for doing this are twofold: it not only provides us a framework to define kinks, but also ensures that any two points that we consider will differ in angle by no more than $ \pi/2 $. This allows to find a lower bound for the probability that the two points belong to the same connected component. If the points in $ V_{x_0,\zeta} $ have angle exceeding $ \pi/2 $, we obtain at most four connected components. 

\begin{definition}[\textbf{Kinks}]\label{D:kinks}
Let $ \Gamma=(V,E) $ be a spoked chain with pole $x_0\notin V$. Suppose that $ x_i,x_j,x_k\in V $ such that $ \hat{x}_i<\hat{x}_j<\hat{x}_k $. We say that $ x_i,x_j $ and $x_k $ form a \textbf{kink} in $ \Gamma $ if the following holds.
\begin{enumerate}
\item $ \widehat{x_ix_jx_k}<\pi/2 $,

\medskip

\item $ \widehat{x_{i^{\sprime}}x_{j^{\sprime}}x_{k^{\sprime}}}\ge \pi/2 $ for all $ x_{i^{\sprime}}, x_{j^{\sprime}}, x_{k^{\sprime}} \in V $ with $ \widehat{x}_{i^{\sprime}}< \widehat{x}_{j^{\sprime}} < \widehat{x}_{k^{\sprime}}  $  and satisfying\\ 
$$ \widehat{x}_i\le \widehat{x}_{i^{\sprime}}<\widehat{x}_{j^{\sprime}}<\widehat{x}_{k^{\sprime}} <\widehat{x}_k \;\mbox{ or }\; \widehat{x}_i<  \hat{x}_{i^{\sprime}}< \hat{x}_{j^{\sprime}} < \widehat{x}_{k^{\sprime}} \le\widehat{x}_k .$$
\end{enumerate}
\end{definition}

\bigskip

\begin{figure}[t]
\centering
\subfigure[Intruding kink]{
\includegraphics[scale=1.2]{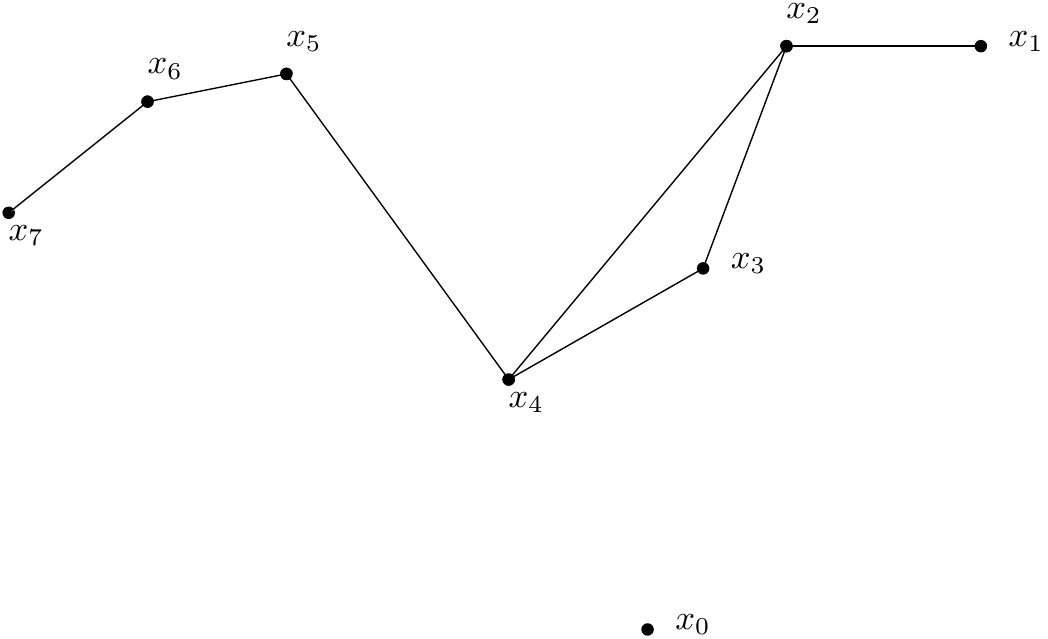}}

\bigskip
\bigskip

\subfigure[Extruding kink]{
\includegraphics[scale=1.0]{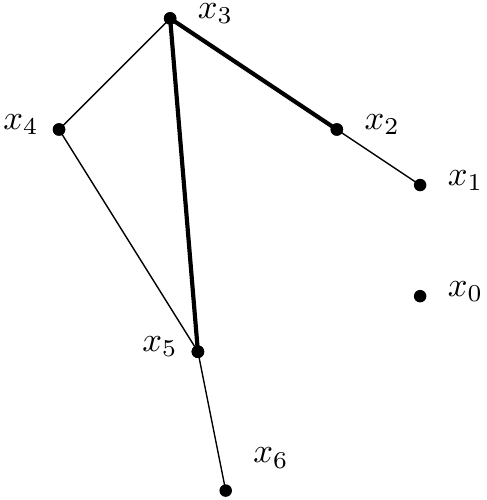}}
\caption{(a) The points $x_2,x_4 $ and $x_5$ form an intruding kink; (b) the points $x_2,x_3 $ and $x_5 $ form an extruding kink}
\label{figc}
\end{figure}

%
%
%
\bigskip

Suppose that $ x_i,x_j $ and $x_k $ form a kink in the spoked chain $ \Gamma=(V,E) $. The kink is called \textbf{intruding} if the line segment $ \overline{x_ix_k} $ lies outside of the induced Polygon $ P(\Gamma,x_0) $ and \textbf{protruding} if it lies inside $ P(\Gamma,x_0) $, see Figure~\ref{figc}. 

\begin{lemma}
Let $ \Gamma=(V,E) $ be a spoked chain with $ V=\{x_1,\ldots,x_n\} $ and  $ \hat{x}_1<\cdots<\hat{x}_n $ and pole $x_0\notin V$. A kink in $ \Gamma $ is either intruding or protruding.
\end{lemma}
\begin{proofsect}{Proof}
Suppose the statement is false. Then there exist $ 1\le i<j<k\le n $ such that $\widehat{x}_i<\widehat{x}_j<\widehat{x}_k $ with $ x_0 $ being the pole such that $ x_i,x_j $ and $x_k $ form a kink in $ \Gamma $ and $ \overline{x_ix_k} $ lies neither inside nor outside of $ P(\Gamma,x_0) $. Let $ U\subset\R^2 $ be the connected component of $ \R^2\setminus\overleftrightarrow{x_ix_k} $ that does not contain $ x_j $. Consider first that $ \overline{x_ix_k} $ lies inside of $ P(\Gamma,x_0) $. It follows that $ \Gamma $ crosses $ \overline{x_ix_k} $ between $ \hat{x}_i $ and $ \hat{x}_k $ and hence  there exists $ x_{j^{\sprime}}\in V\cap U $ with $ \widehat{x}_i<\widehat{x}_{j^{\sprime}}<\widehat{x}_k $. Without loss of generality, let $ \widehat{x}_j<\widehat{x}_{j^{\sprime}}<\widehat{x}_k $. Therefore, as $ x_{j^\prime}\in V\cap U $ one gets $ \widehat{x_ix_jx_{j^{\sprime}}}<\pi/2 $ which contradicts property (2) in Definition~\ref{D:kinks} for the kink formed by $x_i,x_j $ and $x_k $. The second alternative case follows analogously.
\qed
\end{proofsect}

\begin{lemma}\label{L:kink0}
Let $ \Gamma=(V,E) $ be a spoked chain with $ V=\{x_1,\ldots,x_n\} $ and $ \hat{x}_1<\cdots<\hat{x}_n $. 
If $ x_i,x_j $ and $x_k $ form an intruding kink in $ \Gamma $, then
$$
\angle{(\overleftrightarrow{x_ix_{i+1}},\overleftrightarrow{x_{k-1}x_k})}<\pi/2.
$$
\end{lemma}

\begin{proofsect}{Proof}
Let $ x_i,x_j $ and $x_k$ form an intruding kink in $ \Gamma $. Since the kink is intruding, we know that $x_l $ lies in the interior of the triangle $\tau(x_0,x_i,x_k) $ for all $ i<l<k $. Suppose the statement of the Lemma is false, that is,
$$
\angle{(\overleftrightarrow{x_ix_{i+1}},\overleftrightarrow{x_{k-1}x_k})}\ge \pi/2.
$$
This forces either $x_{i+1} $ or $x_{k-1} $ to be in the interior of the triangle $ \tau(x_i,x_j,x_k) $. Without loss of generality, suppose, in fact, that $ x_{i+1} $  is in the interior of that triangle. Therefore, $\widehat{x_{i+1}x_jx_k}<\pi/2 $ which, by Definition~\ref{D:kinks}, contradicts the fact that the points $ x_i,x_j $ and $x_k $ form a kink in $ \Gamma $.
\qed
\end{proofsect}

\begin{figure}[h!]
\caption{Lower bound for the angle $ \widehat{x_0x_{j-1}x_j}$.}\label{angle}
\includegraphics[scale=1.10]{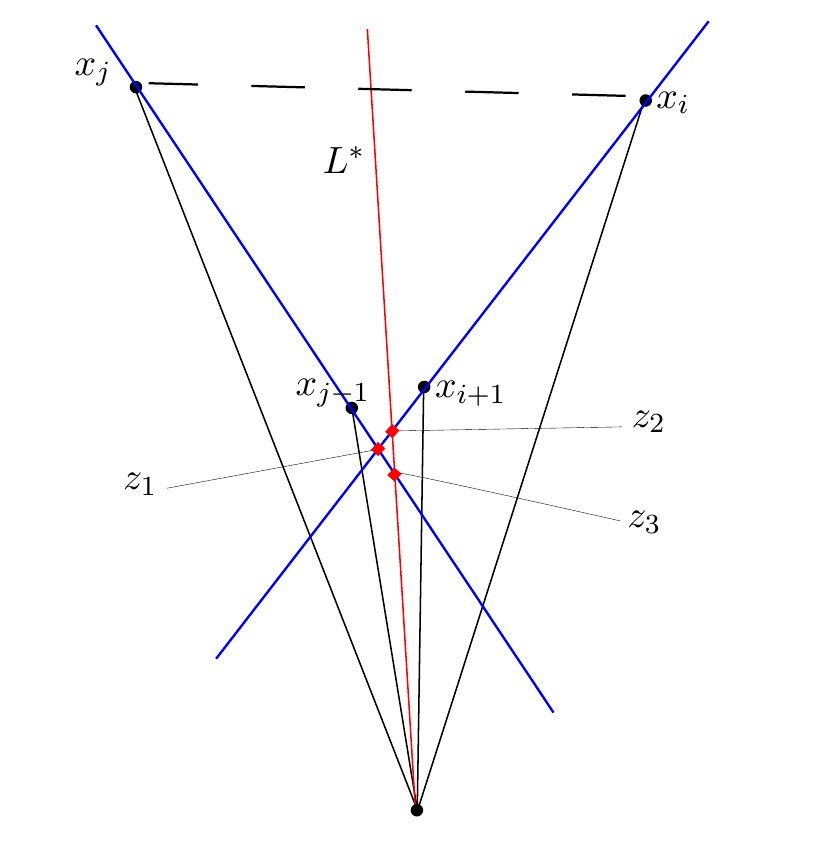}
\end{figure}

\bigskip

\subsection{Intermediary Lemmas}

Kinks of intruding and protruding nature may occur in $ G_B\cap Q_1 $, although their  number  is bounded  from above, see Lemma~\ref{L:kink1} and Lemma~\ref{L:kink2} below.
We shall separate $ G_B\cap Q_1 $ into a finite number of kink-less pieces, each of which a spoked chain, which are easier to work with.  

\begin{lemma}\label{L:kink1}
The number of intruding kinks in $ G_B\cap Q_1 $ is bounded above by $2$.
\end{lemma}

\begin{figure}[h!]
\caption{Lower bound for the angle between kinks of type $2$.}\label{kinkangle}
\includegraphics[scale=0.80]{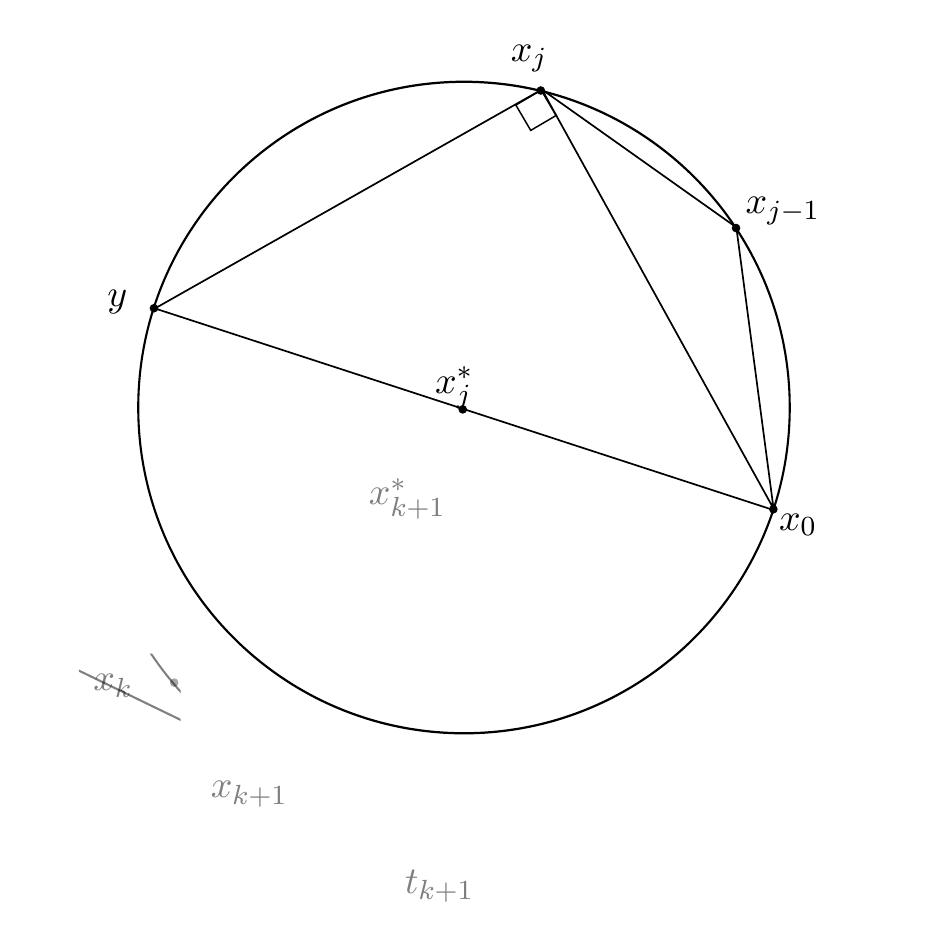}
\end{figure}

\begin{proofsect}{Proof}
We show that the angle between two intruding kinks in a spoked chain is greater than $ \pi/4 $. Since $ G_B\cap Q_1 $ lies in the quadrant $Q_1 $, and is a spoked chain by definition, the statement follows immediately. Let $ \Gamma=(V,E) $ be a spoked chain and order the elements of $V$ such that $ \hat{x}_1<\cdots<\hat{x}_n $. Suppose there is an intruding kink in $ \Gamma $.  By the definition of an intruding kink and Lemma~\ref{L:kink0}, there exist $ 1\le i<n-1 $ and $ i+1<j\le n $, such that
\begin{equation}\label{angleineq}
\angle{(\overleftrightarrow{x_ix_{i+1}},\overleftrightarrow{x_{j-1}x_{j}})}<\pi/2,
\end{equation}
and $ \overline{x_ix_j} $ lies outside of the induced polygon $ P(\Gamma,x_0)$. The straight lines $ \overleftrightarrow{x_ix_{i+1}} $ and $ \overleftrightarrow{x_{j-1}x_j} $ split the plane into four regions. Since the kink is intruding, the point $ x_0 $ must lie in the opposite region to that of the line segment $ \overline{x_ix_j}$. Let $ L^* $ be the radial line of angle $ (\hat{x}_{i+1}-\hat{x}_{j-1})/2 $ and $z_1\in\R^2 $ be the point of the intersection of $ \overleftrightarrow{x_ix_{i+1}} $ and  $ \overleftrightarrow{x_{j-1}x_j} $ and let $z_2,z_3 $ be the points of intersection of $L^* $ with $ \overleftrightarrow{x_ix_{i+1}} $ and $ \overleftrightarrow{x_{j-1}x_j} $ respectively, see Figure~\ref{angle}. Then,
\begin{equation}
\widehat{x_iz_1x_j}+\widehat{x_iz_2x_0}+\widehat{x_0z_3x_j}=2\pi,
\end{equation}
which implies, together with \eqref{angleineq}, that
\begin{equation}
\max\{\widehat{x_iz_2x_0},\widehat{x_0z_3x_j}\}\ge \frac{2\pi-\pi/2}{2}=\frac{3\pi}{4}.
\end{equation}
Without loss of generality, let $\widehat{x_0z_3x_j}\ge 3\pi/4 $. Because $ x_{j-1} $ lies on the line segment $ \overline{z_3x_j} $, it follows that
$$
\widehat{x_0x_{j-1}x_j}\ge\widehat{x_0z_3x_j}\ge 3\pi/4.
$$
Suppose now that there is another intruding kink in $ \Gamma$, formed by the points $x_k,x_l $ and $x_m $ for $ j<k<l<m\le n $. Then, by Lemma~\ref{L:kink0}, we have that
$$
\angle{(\overleftrightarrow{x_kx_{k+1}},\overleftrightarrow{x_{m-1}x_m})}<\pi/2.
$$

\bigskip

\begin{figure}[h!]
\caption{The intruding kink formed by $x_k,x_l $ and $ x_m $.}\label{inkink}
\includegraphics[scale=1.0]{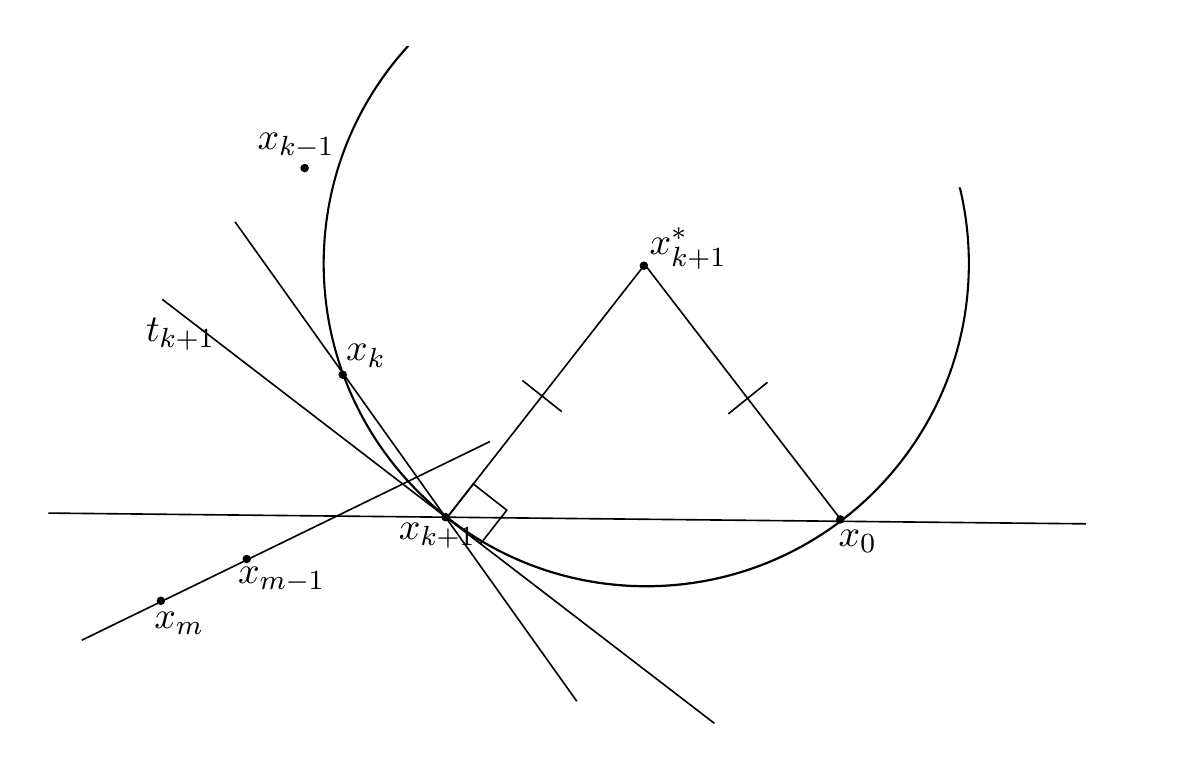}
\end{figure}

\bigskip

Let $ t_{k+1} $ be the tangent to $ \partial B(\tau(x_0,x_k,x_{k+1})) $ at $ x_{k+1} $, see Figure~\ref{inkink}. Then, by noting that 
$$
|V\cap B(\tau(x_0,x_{j-1},x_j))|=0,
$$
which is a consequence of the properties of the Delaunay tessellation (quadratic position), it follows that
\begin{equation}
\begin{aligned}
\angle{(t_{k+1},\overleftrightarrow{x_0x_{k+1}})}&\le \angle{(\overleftrightarrow{x_kx_{k+1}},\overleftrightarrow{x_0x_{k+1}})}\le\angle{(\overleftrightarrow{x_kx_{k+1}},\overleftrightarrow{x_{m-1}x_m})}<\pi/2.
\end{aligned}
\end{equation}
Here, the first inequality follows directly from the definition of the tangent line and the second inequality is a consequence of the fact that $ \hat{x}_{k+1}<\hat{x}_{m-1}<\hat{x}_m $. For $ 1\le r<n $, let $ x^*_{r+1} $ denote the centre of the circumcircle of the triangle $ \tau(x_0,x_r,x_{r+1})\in\Del_3(V\cup\{x_0\}) $. Since the points $ \{x_0,x^*_{k+1},x_{k+1}\} $ form an isosceles, see Figure~\ref{inkink}, we can conclude that
\begin{equation}\label{anglek}
\hat{x}_{k+1}-\hat{x}^*_{k+1}=\pi/2-\angle{(t_{k+1},\overleftrightarrow{x_0x_{k+1}})}>0.
\end{equation}
Let $ y$ be the antipodal point to $x_0 $ on the circumscribed ball of the triangle $ \tau(x_0,x_{j-1},x_j) $ in $ \R^2 $. Since, $|x_0-y| $ is equal to the diameter of that circle, it follows that $ \widehat{x_0x_jy}=\pi/2 $, see Figure~\ref{kinkangle}. The points $x_0,x_{j-1},x_j $ and $y$ form a cyclic quadrilateral. Using 
$ \widehat{x_0x_{j-1}x_j}\ge\widehat{x_0z_3x_j}\ge 3\pi/4 $ from above, and the fact that opposite angles of a cyclic quadrilateral add up to $ \pi $, we see that $ \widehat{x_0yx_j}\le \pi/4$.
Hence, by \eqref{anglek}
$$
\hat{x}_{k+1}-\hat{x}_j>\hat{x}^*_{k+1}-\hat{x}_j\ge\hat{x}^*_j-\hat{x}_j=\widehat{yx_0x_j}=\pi-\widehat{x_0x_jy}-\widehat{x_0yx_j}\ge \pi/4,
$$ 
where the second inequality is due to a further property of the Delaunay structure, see Lemma~\ref{L:E} in appendix~\ref{appE}. This implies that the angle between intruding kinks must be greater than $ \pi/4 $.
\qed
\end{proofsect}

\begin{lemma}\label{L:kink2}
There are no protruding kinks in $ G_B\cap Q_1 $.
\end{lemma}
\begin{proofsect}{Proof}
We order the elements of $ V_B=\{x_1,\ldots,x_n\} $ such that $ \hat{x}_1<\cdots<\hat{x}_n $. Suppose we have a protruding kink with pole $x_0\notin V$, then we have for some $ 1\le i<j<k\le n $,
$$
\widehat{x_ix_jx_k}<\pi/2.
$$
The pair $\{x_i,x_k\} $ does not form an edge of $E_B $, therefore, by the properties of the Delaunay tessellation, $x_j$ lies inside the circumcircle $\partial B(\tau(x_0,x_i,x_k)) $ of the triangle $ \tau (x_0,x_i,x_k) $. The line segment $ \overline{x_ix_k} $ is a chord which splits the ball $ B(\tau(x_0,x_i,x_k))$  into two regions. Since we have a protruding kink, $\overline{x_ix_k} $ lies inside the induced polygon $ P(G_B\cap Q_1;x_0) $, and so $x_j$ does not lie in the same region as $x_0 $. The angle $ \widehat{x_ix_jx_k} =\pi/2 $ once the point $x_j $ lies on the boundary of the ball $ B(\tau(x_0,x_i,x_k)) $ for the case $  \widehat{x_ix_0x_k}=\pi $, and due to the fact that  $ \widehat{x_ix_jx_k}<\pi/2 $ we get that $ \widehat{x_ix_0x_k}\ge \pi/2 $, and hence, there are no protruding kinks in $ G_B\cap Q_1 $.
\qed
\end{proofsect}

\bigskip

\subsection{Edge drawing}\label{S:edge}
Before we can finish the proof of Theorem~\ref{THMc} in Section~\ref{S:final} below we need two more results to gain some control over the edge drawing mechanism in $ E^{\ssup{\rm ext}}_{x_0,\zeta} $. Denote by $ \widetilde{\nu}_{\zeta} $ the edge drawing mechanism with probability
\begin{equation}
\widetilde{p}_2(\eta_{xy})=\begin{cases}\frac{\1\{|x-y|\le R\}}{\frac{q}{\beta}|x-y|^4+1}  & \mbox{ if } \eta_{xy}\in\Del_2(\zeta)\cap (E_{\R^2}\setminus E_{\L^{\rm c}}),\\
\1_{\Del_2(\zeta)}(\eta_{xy}) & \mbox{ if } \eta_{xy}\in E_{\L^{\rm c}}\,,\end{cases} 
\end{equation}
and denote $ \widetilde{\nu}^{\ssup{\rm ext}}_{\zeta} $ the corresponding edge drawing mechanism on $ E^{\ssup{\rm ext}}_{x_0,\zeta} $.

\begin{lemma}
For all $ q\ge 1 $, $ \L\Subset\R^2 $  and all $ \zeta\in\O_{\L,\omega} $, 
$$
\mu^{\ssup{q}}_{\L,\zeta}   \succcurlyeq \widetilde{\nu}_{\zeta}.
$$
\end{lemma}

\begin{proofsect}{Proof}
It suffices to show that for all edges $ \eta\in\Del_2(\zeta) $ with $ \eta\notin E_{\L^{\rm c}}$,
\begin{equation}\label{ineqfinal}
\frac{p(\eta)}{q(1-p(\eta))}\ge\frac{\widetilde{p}_2(\eta)}{(1-\widetilde{p}_2(\eta))}.
\end{equation}
Recall from \eqref{tiledrawing} that
$$
p(\eta)=\frac{1}{\frac{1}{\beta}\ell(\eta)^4+1}.
$$
Thus, if $ \ell(\eta)>R $, then $ \widetilde{p}(\eta)=0 $ and the inequality \eqref{ineqfinal} trivially holds. Suppose that $ \ell(\eta)<R $, then, in fact, we also have 
$$
\frac{p(\eta)}{q(1-p(\eta))}=\frac{\widetilde{p}_2(\eta)}{(1-\widetilde{p}_2(\eta))}.
$$
Henceforth, \eqref{ineqfinal} holds for all $ \eta\in\Del_2(\zeta) $.
\qed
\end{proofsect}

Note that $E_B $ is not necessarily a subset of $ E^{\ssup{\rm ext}}_{x_0,\zeta} $, in fact, they belong to different Delaunay tessellations
$$
E_B\subset\Del_2(V_B\cup\{x_0\}) \,\mbox{ and }\, E^{\ssup{\rm ext}}_{x_0,\zeta}\subset\Del_2(\zeta\cup\{x_0\}).
$$
We therefore introduce another edge drawing mechanism, but this time on $E_B $. Let $ \mu^* $ denote the distribution of the random edge configurations $ \{\eta\in E_B\colon \upsilon(\eta)=1\} $,
where $ ((\upsilon(\eta))_{\eta\in E_B} $ are independent Bernoulli random variables with probability
\begin{equation}\label{pstar}
\P(\upsilon(\eta)=1)=p^*(\eta)=\frac{\1\{\ell(\eta)\le \frac{2}{\pi}\wedge R\}\1\{|\hat{x}-\hat{y}|\le\frac{\pi}{2}\}}{\frac{q}{\beta}\big(\frac{\pi}{2}\ell(\eta)^4\big)^4+1}\1_{E_B}(\eta),\quad\mbox{ for } \eta=\eta_{xy}.\end{equation}
We now compare the probability that two points are connected  with respect to $ \widetilde{\nu}_\zeta $ and with respect to $ \mu^* $.

\begin{lemma}
Pick $ \zeta\in\O_{\L,\omega} $. Let $ \eta_{xy}\in E_B $ and let $ x\leftrightarrow y $ denote the event that $x $ and $y $ lie in the same connected component of $ (\zeta,E)$, where $E$ is a $\widetilde{p}_2$-thinning of the edge set $ E^{\ssup{\rm ext}}_{x_0,\zeta} $. Then,
\begin{equation}\label{thinning}
\widetilde{\nu}^{\ssup{\rm ext}}_{\zeta} (x\leftrightarrow y)\ge p^*(\eta_{xy}).
\end{equation}
\end{lemma}

\begin{proofsect}{Proof}
By the definition of $p^*$, \eqref{thinning} follows for $x,y\in V_B $ with $ |x-y|>\frac{2}{\pi}\wedge R $ or with $|\hat{x}-\hat{y}|>\frac{\pi}{2} $. Therefore, we assume that $ |x-y|\le\frac{2}{\pi} \wedge R $ and $ |\hat{x}-\hat{y}|\le \frac{\pi}{2} $.

\noindent \textbf{Case I:} If $ \eta_{xy}\in E^{\ssup{\rm ext}}_{x_0,\zeta} $, we get
$$
\widetilde{\nu}^{\ssup{\rm ext}}_\zeta(x\leftrightarrow y)\ge \widetilde{p}_2(\eta_{xy})=\frac{1}{\frac{q}{\beta}\ell(\eta_{xy})^4+1}\ge \frac{1}{\frac{q}{\beta}(\frac{\pi}{2}\ell(\eta_{xy}))^4+1} =  p^*(\eta_{xy}).
$$

\noindent\textbf{Case II:} If $ \eta_{xy}\not\in E^{\ssup{\rm ext}}_{x_0,\zeta} $, the proof is no longer straightforward and will take some care. Since $ \eta_{xy}\not\in E^{\ssup{\rm ext}}_{x_0,\zeta} $, and $x,y\in V_B $ there exists $z\in \zeta \cap V_B^{\rm c} $, such that $ \eta_{zx_0}\in\Del_2(\zeta\cup\{x_0\})$. This implies that $ z\in V_{x_0,\zeta}\setminus V_B $ and $ \hat{x}<\hat{z}<\hat{y} $. We now check whether $ \eta_{xz},\eta_{zy}\in E^{\ssup{\rm ext}}_{x_0,\zeta} $. If they are not, we find more points  $ z\in V_{x_0,\zeta}\setminus V_B $ with $ \hat{x}<\hat{z}<\hat{y} $. Therefore, there exists a (finite) sequence  $z_1,\ldots, z_n\in V_{x_0,\zeta}\setminus V_B $ with $ \hat{x}<\hat{z}_1<\cdots<\hat{z}_n<\hat{y} $ such that
$$
\eta_{xz_1},\eta_{z_1,z_2},\ldots,\eta_{z_ny}\in E^{\ssup{\rm ext}}_{x_0,\zeta}.
$$

The event that each of these edges is open implies the event that $ x$ and $y$ belong to the same connected component of open edges, hence
\begin{equation}
\widetilde{\nu}^{\ssup{\rm ext}}_\zeta(x\leftrightarrow y)\ge \widetilde{p}_2(\eta_{xz_1}) \widetilde{p}_2(\eta_{z1z_2})\cdots \widetilde{p}_2(\eta_{z_ny}).
\end{equation}
For any two points $x_1,x_2\in\zeta $ with $ \hat{x}_1<\hat{x}_2 $, define $ C^{x_0}_{x_1,x_2} $ to be the arc on the circumcircle $ \partial B(\tau(x_1,x_2,x_0)) $ of the triangle $ \tau(x_1,x_2,x_0) $ between $ x_1 $ and $x_2 $,  and define $U_{x_1,x_2} $ to be the subset of $ \R^2 $ bounded by this arc $ C^{x_0}_{x_1,x_2} $ and $ \overline{x_1x_2} $, that is, the convex hull of $ C^{x_0}_{x_1,x_2} $. Let $n=\#\{z\in V_{x_0,\zeta}\colon \hat{x}<\hat{z}<\hat{y}\} $. We claim that
\begin{equation}\label{claim}
L(C^{x_0}_{x,z_1})+\cdots+L(C^{x_0}_{z_n,y})\le L(C^{x_0}_{x,y}),\quad n\in\N,
\end{equation}
where $ L(C) $ denotes the length of the arc $C$. We will prove the claim \eqref{claim} below after we finish the proof of the statement in the lemma. 
By our assumption that  $ |x-y|\le\frac{2}{\pi} \wedge R $ and $ |\hat{x}-\hat{y}|\le \frac{\pi}{2} $, it follows that $ L(C^{x_0}_{x,y})\le 1 $. To see this note that with $r $ being the radius of the circumcircle and with $ \theta=\widehat{xx_0y}=|\hat{x}-\hat{y}| $,
$$
L(C^{x_0}_{x,y})=2r\theta=2\theta\frac{|x-y||x_0-x||x_0-y|}{4\,\area(\tau(x,x_0,y))}=\frac{|x-y|\theta}{\sin(\theta)}\le 1,
$$
where we used that $ \sin(\theta)/\theta >2/\pi $ for $ \theta\in (0,\pi/2)$.  With our claim  \eqref{claim} we obtain 
$$
L(C^{x_0}_{x,z_1})+\cdots+L(C^{x_0}_{z_n,y})\le 1.
$$
Obviously, this shows that
\begin{equation}\label{length}
|x-z_1|+|z_1-z_2|+\cdots +|z_{n-1}-z_n|+|z_n-y|\le 1. 
\end{equation}
Now choose $ \beta >0 $ such that $ q/\beta <1 $. Then for $ a,b\in\R $ with $ 0\le a\le b\le 1 $, we have the following simple fact
\begin{equation}\label{fact}
\Big(\frac{1}{\frac{q}{\beta}a^4+1}\Big)\Big(\frac{1}{\frac{q}{\beta}b^4+1}\Big)=\frac{1}{\frac{q}{\beta}\big(\frac{q}{\beta}a^4b^4+a^4+b^4\big)+1}\ge\frac{1}{\frac{q}{\beta}(a+b)^4+1},
\end{equation}
where the inequality follows because of $ q/\beta<1 $ and the given constraints on $a $ and $b$. Hence, using \eqref{length}, we obtain
\begin{equation}
\begin{aligned}
\widetilde{p}_2(\eta_{xz_1}) &\widetilde{p}_2(\eta_{z1z_2})\cdots \widetilde{p}_2(\eta_{z_ny})\ge \Big(\frac{1}{\frac{q}{\beta}|x-z_1|^4+1}\Big)\cdots \Big(\frac{1}{\frac{q}{\beta}|z_n-y|^4+1}\Big)\\
&\ge \frac{1}{\frac{q}{\beta}\big(|x-z_1|+|z_1-z_2|+\cdots +|z_{n-1}-z_n|+|z_n-y|\big)^4+1}\\
&\ge \frac{1}{\frac{q}{\beta}L(C^{x_0}_{x,y})^4+1}\ge \frac{1}{\frac{q}{\beta}\big(\frac{\pi}{2}|x-y|\big)^4+1}=p^*(\eta_{xy}),
\end{aligned}
\end{equation}
where the second inequality results from repeated use of relation \eqref{fact} with $a=|z_i-z_{i+1}| $ and $b=|z_j-z_{j+1}| $. We conclude with the statement in the lemma. 

\noindent We are left to verify the claim \eqref{claim}: Suppose there exists $z\in V_{x_0,\zeta}\setminus V_B $ such that $ \eta_{xz},\eta_{zy}\in E^{\ssup{\rm ext}}_{x_0,\zeta} $. Since $ z\not\in B $, it must be true that $z\in U_{xy} $. Therefore, by a direct application of Theorem~\ref{T:E}, we have
$$
L(C^{x_0}_{x,z})+L(C^{x_0}_{z,y})\le L(C^{x_0}_{x,y}),
$$
and the claim follows for $n=1 $. We shall proceed by induction with respect to $n\in\N $. Assume the claim holds for $ n=k-1 $. There exist $ z_1,\ldots,z_k\in V_{x_0,\zeta}\setminus V_B $ such that $ \hat{x}_1<\cdots<\hat{x}_k $ and $ \eta_{x,z_1},\ldots,\eta_{z_ky}\in E^{\ssup{\rm ext}}_{x_0,\zeta} $. Let
$$
i=\argmax_{1\le j\le k} |z_j-\overline{xy}|.
$$
It follows that $z_i\in U_{z_{i-1},z_{i+1}} $,  where, for convenience, we write $ z_0=x $ and $z_{k+1}=y $. By Theorem D.2 again,
\begin{equation}
L(C^{x_0}_{z_{i-1},z_i})+L(C^{x_0}_{z_i,z_{i+1}})\le L(C^{x_0}_{z_{i-1},z_{i+1}}).
\end{equation}
By  changing the notation  $ z^{\sprime}_j=z_j $ for $ 1\le j <i $ and $ z^{\sprime}_j=z_{j+1} $ for $ i\le j\le k-1 $, it follows from the previous inequality that
$$
L(C^{x_0}_{x,z_1})+\cdots +L(C^{x_0}_{z_k,y})\le L(C^{x_0}_{xz^{\sprime}_a})+\cdots +L(C^{x_0}_{z^{\sprime}_{k-1},y}),
$$
and hence, by our assumption for $ n=k-1 $, we conclude with the statement of the claim \eqref{claim}.
\qed
\end{proofsect}

\bigskip

\begin{figure}[h!]
\caption{The sectors $S_i $ of a spoked chain in $Q_1 $.}\label{spokedsector}
\includegraphics[scale=1.08]{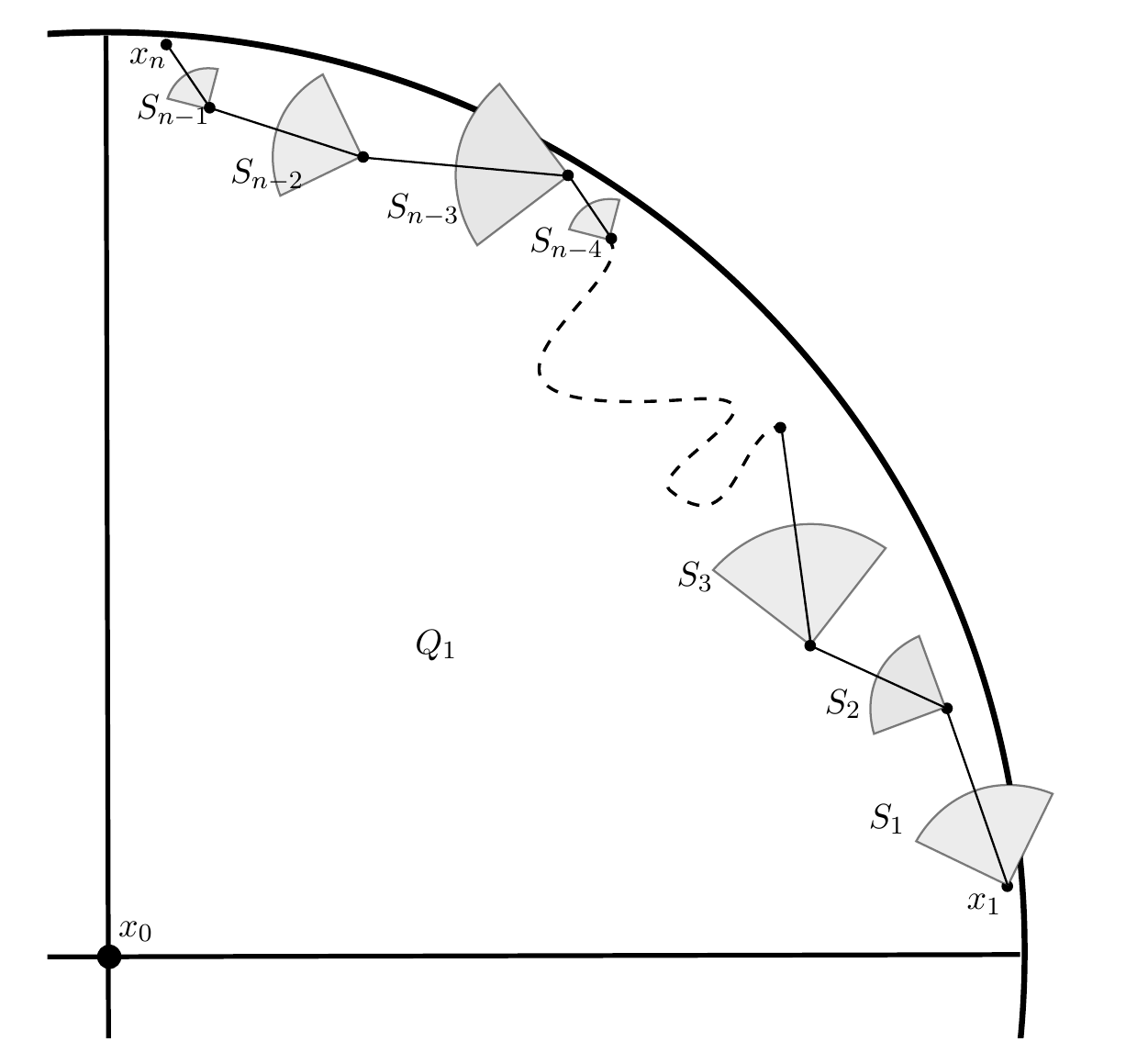}
\end{figure}

\bigskip

\begin{lemma}\label{L:deltabound}
Let $\delta>0 $ and $ \Gamma=(V,E) $ be a spoked chain with $V\subset Q_1 $. If $ \Gamma $ does not contain a kink, then the number of edges in $E$ with length greater than $ 2\delta $ is at most $ 6(\frac{R}{\delta})^2$. \end{lemma}

\begin{proofsect}{Proof}
Let $ V=\{x_1,\ldots,x_n\} $ with $ \hat{x}_1<\cdots<\hat{x}_n $ be given. For $1\le i <n $, let $D_i\subset\R^2 $ be the disc of radius $|x_i-x_{i+1}|/2 $ centred at $x_i $. Let $S_i\subset\R^2 $ be the sector of $D_i $ with interior angle $\pi/2 $ and line of symmetry $ \overline{x_ix_{i+1}} $. We claim the following,
\begin{equation}
S_i\cap S_j=\emptyset\,\mbox{ for } i\not= j, 1\le i,j\le n,\label{claimone}
\end{equation} and 
\begin{equation}
\bigcup_{i=1}^{n-1}S_i  \subset \Ucal_{\frac{R}{\sqrt{2}}}(Q_1),\label{claim2}
\end{equation}
where $ \Ucal_{\frac{R}{\sqrt{2}}}(Q_1)=\{x\in\R^2\colon \dist(x,Q_1)\le \frac{R}{\sqrt{2}}\} $ is the $ R/\sqrt{2} $ neighbourhood of the sector $Q_1 $.
Assume our claim \eqref{claim2} is true, the sum of the areas of the sectors $S_i $ must not exceed the area of $\Ucal_{\frac{R}{\sqrt{2}}}(Q_1) $ which is less than $ \frac{3}{2}\pi R^2 $, see Figure~\ref{spokedsector}. Now each edge $\eta\in E $ of length greater than $ 2\delta $ contributes a sector of area greater than $ \pi/4 \delta^2 $, therefore, the maximum number of such edges in $ \Gamma $ is simply
$$
\frac{(3/2)\pi R^2}{\pi/4 \delta^2}=6\Big(\frac{R}{\delta}\Big)^2,
$$
which gives the result. We are left to prove our claim \eqref{claim2} above. Pick $ x_i\in V $ and let $\ell_1 $ be the image of the line $\overleftrightarrow{x_ix_{i+1}} $ under a rotation with an angle $\pi/2 $, centred at point $ x_{i+1} $. There are exactly two connected components of $\R^2\setminus\ell_1 $. Let $ U $ denote the component that contains $ x_i $. Now suppose $x_k\in U $ for some $i+1<k\le n $. This implies that $ \widehat{x_ix_{i+1}x_k}<\pi/2 $. Then, by Definition~\ref{D:kinks}, this contradicts the fact that $ \Gamma $ does not contain a kink. Therefore, $ x_k\in U^{\rm c} $ for all $i+1<k\le n $. Let $\ell_2 $ and $ \ell_3 $ be the images of the half line $ \overleftarrow{x_ix_{i+1}} $ under rotations, centred at $x_{i+1} $, of angles $\pi/4 $ and $ -\pi/4 $ respectively, see Figure~\ref{sector}.
\bigskip

\begin{figure}[h!]
\caption{The point $x^{\sprime}_k $ is the first time after $ x_{i+1} $ that the chain enters $U$.}\label{sector}
\includegraphics[scale=0.90]{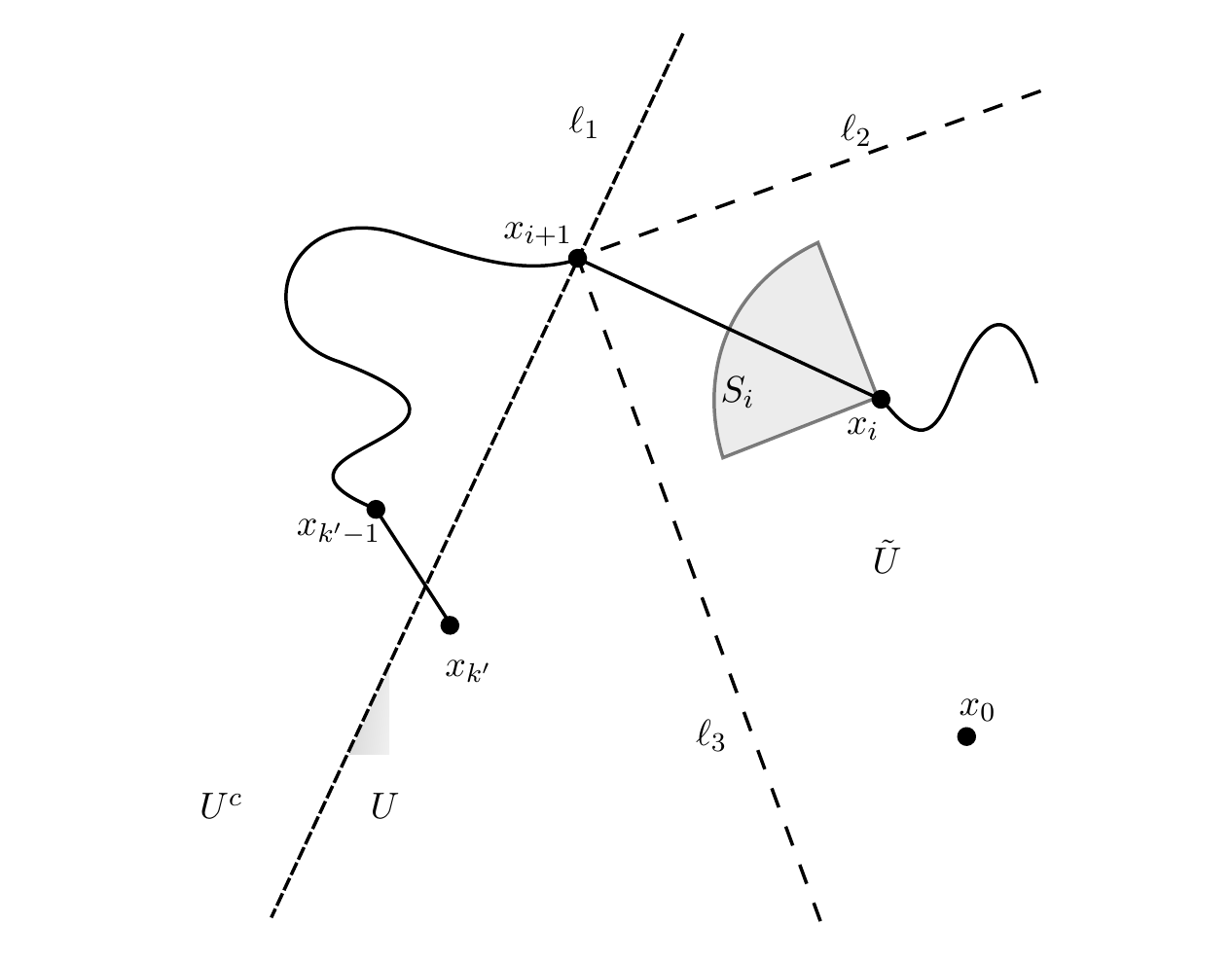}
\end{figure}

\bigskip

Again, there are two connected components of $ \R^2\setminus(\ell_2\cup\ell_3) $. Let $ \widetilde{U} $ denote the one that contains $x_i$. Now \eqref{claimone}  follows by noticing that $S_i\subset\widetilde{U} $ and $S_k\subset\widetilde{U}^{\rm c} $ for all $ i+1<k\le n $. Claim  \eqref{claim2} follows easily  $ S_i\subset D_i $ for all $ 1\le i<n $ and that the maximal radius for $D_i $ is half the maximal edge length, which is $ \sqrt{2}R $.

 \qed
\end{proofsect}
\subsection{Final proof of Theorem~\ref{THMc}}\label{S:final}
Recall that we split $B=B_R(x_0) \subset\R^2 $ into four quadrants, $Q_i\subset\R^2,i=1,2,3,4$. Now $ G_B\cap Q_1 $ contains all vertices and edges of $G_B $ that lie wholly in $Q_1 $. By construction, $ G_B\cap Q_1 $ is a spoked chain. It follows from Lemma~\ref{L:kink1} and Lemma~\ref{L:kink2} that there are at most $2$ intruding kinks in the spoked chain $ G_B\cap Q_1 $ and no single protruding kink. For each intruding kink $x_i,x_j,x_k $, we remove the edge $ \eta_{x_jx_{j+1}} $ from $ G_B\cap Q_1 $. Since removing an edge anywhere except from the end of the spoked chain will result in leaving two spoked chains, we are left with at most $3 $ spoked chains in $Q_1 $. Importantly, none of these contain an intruding or protruding kink. Let denote  $ \Gamma=(V^\Gamma,E^\Gamma) $ one of these kink-less spoked chains in $Q_1 $. We denote (compare with Theorem~\ref{THMc}) $ N^{\ssup{\rm cc}}_{\Gamma}(\zeta,E) $ to be the number of connected components (clusters) of $ (\zeta,E) $ that intersect $V^\Gamma $. We endeavour to bound the expectation of $ N^{\ssup{\rm cc}}_{x_0}(\zeta,\cdot)$ with respect to the edge drawing mechanism $ \mu^{\ssup{q}}_{\ssup{\rm ext},\zeta} $ on $ E^{\ssup{\rm ext}}_{x_0,\zeta} $ given in \eqref{edgem}.  To conclude the Theorem, we shall use 
$$
\int\; N^{\ssup{\rm cc}}_{x_0}(\zeta,E)\,\mu^{\ssup{q}}_{\ssup{\rm ext},\zeta}(\d E)\le 12\int\:N_\Gamma^{\ssup{\rm cc}}(\zeta,E)\, \mu^{\ssup{q}}_{\ssup{\rm ext},\zeta}(\d E),
$$
where the factor $12$ is considering at most three kinkless spoked chains in each of the four quadrants. Order the elements in $V_\Gamma=\{x_1,\ldots,x_n\} $ such that  $ \hat{x}_1<\cdots<\hat{x}_n $. Recall that $ \{x\leftrightarrow y\} $ denotes the event that $x$ and $y $ belong to the same cluster of $ (\zeta,E) $ and notice that 
\begin{equation}\label{int}
\begin{aligned}
\int\:N_\Gamma^{\ssup{\rm cc}}(\zeta,E)\, \mu^{\ssup{q}}_{\ssup{\rm ext},\zeta}(\d E)& \le 1+\sum_{j=1}^{n-1}\big(1-\mu^{\ssup{q}}_{\ssup{\rm ext},\zeta}(\{x_j\leftrightarrow x_{j+1}\})\big)\\&
 \le 1+\sum_{j=1}^{n-1}\big(1-\widetilde{\nu}^{\ssup{\rm ext}}_\zeta(\{x_j\leftrightarrow x_{j+1}\})\big)\\
&\le 1+ \sum_{j=1}^{n-1}\big(1-p^*(\eta_{x_jx_{j+1}})\big)\\
&\le 1+  \sum_{\eta\in E^\Gamma}\big(1-p^*(\eta)\big).
\end{aligned}
\end{equation}

We partition the edge set $E^\Gamma $ of the spoked chain $ \Gamma $ into subsets of edges according to their lengths. Let
$$
\begin{aligned}
E_1&=\{\eta_{xy}\in E^\Gamma\colon |x-y|>\frac{2}{\pi}\wedge R\},\\[1ex]
E_i&=\Big\{\eta_{xy}\in E^\Gamma\colon \frac{\frac{2}{\pi}\wedge R}{i}<|x-y|\le\frac{\frac{2}{\pi}\wedge R}{i-1}\Big\},\quad i\ge 2,i\in\N.
\end{aligned}
$$
By recalling that
$$
p^*(\eta_{xy})=\frac{\1\{\ell(\eta)\le \frac{2}{\pi}\wedge R\}\1\{|\hat{x}-\hat{y}|\le\frac{\pi}{2}\}}{\frac{q}{\beta}\big(\frac{\pi}{2}\ell(\eta)^4\big)^4+1}\1_{E_B}(\eta)
$$ 
from \eqref{pstar}, we see that $1-p^*(\eta)=1 $ for all $ \eta\in E_1 $. However, since $ \Gamma $ is contained in $Q_1 $, and henceforth $ |\hat{x}-\hat{y}|<\frac{\pi}{2} $, we have
$$
1-p^*(\eta_{xy})=\frac{1}{\frac{q}{\beta}\big(\frac{\pi}{2}|x-y|\big)^4+1}
$$ for all $ \eta_{xy}\in E_i, i\ge 2 $. Let $r:=1\wedge \frac{R\pi}{2} $. Then, considering $ \frac{2r}{\pi i}< |x-y|\le \frac{2r}{\pi(i-1)} $ for all $ \eta_{xy}\in E_i $ and noticing that $ \bigcup_{i=1}^\infty E_i=E^\Gamma $ and $ E_i\cap E_j=\emptyset $ for $ i\not= j $, it follows readily that
$$
\begin{aligned}
\sum_{\eta\in E^\Gamma}\big(1-p^*(\eta)\big)&=\sum_{i=1}^\infty\sum_{\eta\in E_i} \big(1-p^*(\eta)\big)\le \sum_{\eta\in E_1} 1+\sum_{i=2}^\infty\sum_{\eta\in E_i}\Big(\frac{1}{\frac{\beta}{q}\big(\frac{i-1}{r}\big)^4+1}\Big)\\
&\le 6R^2\pi^2r^{-2}+\sum_{i=2}^\infty\big(6R^2\pi^2i^2r^{-2}\big) \Big(\frac{1}{\frac{\beta}{q}\big(\frac{i-1}{r}\big)^4+1}\Big)\\
&\le 6R^2\pi^2r^{-2}\Big(1+\frac{qr^2}{\beta}\sum_{i=2}^\infty\frac{i^2}{(i-1)^4}\Big), 
\end{aligned}
$$
where the second inequality comes from an application of Lemma~\ref{L:deltabound}. We use that
\begin{equation}
\sum_{i=2}^\infty\frac{i^2}{(i-1)^4}\le \sum_{i=2}^\infty\frac{4}{(i-1)^2}=4\sum_{i=1}^\infty\frac{1}{i^2}=\frac{2}{3}\pi^2.
\end{equation}
Combining all our previous steps we obtain from \eqref{int}  that
$$
\begin{aligned}
\int\:N_\Gamma^{\ssup{\rm cc}}(\zeta,E)\, \mu^{\ssup{q}}_{\ssup{\rm ext},\zeta}(\d E) \le 1+6R^2\pi^2r^{-2}\big(1+\frac{2q\pi^2r^2}{3\beta}\big)\\
\end{aligned}
$$
We finish the proof of Theorem~\ref{THMc} by setting $ \alpha=\alpha(R,q,\beta)= 1+6R^2\pi^2r^{-2}\big(1+\frac{2q\pi^2r^2}{3\beta}\big)$. Note that for given $ \beta $ and $ q $ the function $ \alpha(R,q,\beta) $ grows quadratic in the finite range radius $R$ (in both cases $ R\pi/2> 1 $ and $ R\pi/2< 1 $). Furthermore, 
$$
\lim_{\beta\to\infty}\alpha(R,q,\beta)=1+6R^2\pi^2r^{-2}.
$$

\qed

\section{Proofs}\label{proofs} 
This section delivers the remaining open proofs of our results. We first establish the existence of  Gibbs measures. In Section~\ref{finish} we finally finish the proof of Theorem~\ref{THM-mainedge}. 
\subsection{Existence of Gibbs measures}
To show the existence of Gibbs measures (Proposition~\ref{Propexist}) for our Delaunay Potts model we follow \cite{DDG12}.   The potential $ \phi_\beta $  depends solely on the individual  Delaunay hyperedges in $ \Del_2(\bzeta) $, of a marked configuration $ \bzeta $. Every marked hyperedge $ \beeta\in\Del_2(\bzeta) $ has the so-called finite horizon $\overline{B}(\eta,\zeta) $,  where $ B(\eta,\zeta)$ is the open ball with $ \partial B(\eta,\zeta)\cap\zeta=\eta $ that contains no points of $ \zeta$.  Thus $ \phi_\beta $ satisfies  the range condition (R) in \cite{DDG12}, see \cite[Proposition~4.1 \& 4.3]{DDG12}, with finite horizon being the ball $ \overline{B}(\eta,\zeta) $.  The finite-horizon property of a general hyperedge potential $ \varphi\colon\bO\times\bO\to\R  $ says that for each pair $ (\beeta,\bzeta) $ with $ \beeta\in\Del_2(\bzeta) $ there  exists some $ \Delta\Subset\R^2 $ such that for the pair $ (\beeta,\widetilde{\bzeta}) $ with $ \beeta\in\Del_2(\widetilde{\bzeta}) $ we have that $ \varphi(\beeta,\bzeta)=\varphi(\beeta,\widetilde{\bzeta}) $ when $ \widetilde{\bzeta}=\bzeta $ on $ \Delta\equiv\overline{B}(\eta,\zeta) $.   

The second requirement for existence of Gibbs measures is the  stability condition (S). A hyperedge potential  is called stable if there is a lower bound for the Hamiltonian for any $ \L\Subset\R^2 $, and, as $ \phi_\beta(\ell)\ge 0 $ for all $ \ell\ge 0 $,  the stability condition (S) is satisfied. The third condition to be checked  is a partial complementary upper bound for the Hamiltonian in any $ \L\Subset \R^2 $. This is a bit more involved,  and we shall first define appropriate configurations, the so-called \textit{pseudo-periodic marked configurations}. We consider the partition of $ \R^2 $ as given in Appendix~\ref{pseudoperiodic}. Note that in Appendix~\ref{pseudoperiodic} we have introduced a length scale $ \ell >0 $ which is not necessary for the existence proof as we can put $ \ell=1 $. We let $ B(0,r) $ be an open ball of radius $ r\le \rho_0\ell $,   where we choose $ \rho_0\in (0,1/2) $ sufficiently small such that  $ B(0,r)\subset\Delta_{0,0} $. Note that 
$$ B^r:=\{\zeta\in\O_{\Delta_{0,0}}\colon \zeta=\{x\}\mbox{ for some } x\in B(0,r)\} $$ is a measurable set of $\O_{\Delta_{0,0}}\setminus\{\emptyset\} $. Then
$$
\Gamma^r=\{\om\in\O\colon\theta_{Mz}(\om_{\Delta_{k,l}})\in B^r\mbox{ for all } (k,l)\in\Z^2\} 
$$ is the set of pseudo-periodic configurations \eqref{ppc}. These configurations are not marked yet. The reason is that when a point is shifted its mark remains unchanged. Thus we define the set of pseudo-periodic marked configurations as 
$$
\boldsymbol{\Gamma}^r=\{\bo=(\om^{\ssup{1}},\ldots,\om^{\ssup{q}})\colon \om^{\ssup{i}}\in\Gamma^r\mbox{ for all }i\in M_q\}.
$$

The required control of the Hamiltonian from above will be achieved by the following properties. As our hyperedge potential depends only on the single hyperedge the so-called uniform confinement (see \cite{DDG12}) is trivially satisfied. In addition,  we need the uniform summability, that is,
$$
c_r:=\sup_{\bzeta\in\boldsymbol{\Gamma}^r}\sum_{\beeta\in\Del_2(\bzeta)\colon \eta\cap\Delta\not=\emptyset}\frac{\phi_\beta(\ell(\eta))(1-\delta_\sigma(\beeta))}{\#\widehat{\eta}}<\infty,
$$
where $ \widehat{\eta}=\{(k,l)\in\Z^2\colon \eta\cap \Delta_{k,l}\not=\emptyset\} $ and where $ \Delta=\Delta_{0,0}$. The length $ \ell(\eta) $ of any $ \eta\in\Del_2(\zeta)\cap \Delta $ when $ \zeta $ is any pseudo-periodic configuration  satisfies  $ \ell(1-2\rho_0)\le \ell(\eta)\le \ell(1+2\rho_0) $.  There are at most six edges from the centre ball in $ \Delta=\Delta_{0,0} $ and each Delaunay edge touches exactly two cells and thus $ \widehat{\eta}=2 $. We obtain an upper bound for each edge by considering the shortest possible length for each edge, that is,
$$
c_r=3\log\Big(\frac{(\ell(1-2\rho_0))^4+\beta}{(\ell(1-2\rho_0))^4}\Big)<\infty.
$$
We need furthermore the so-called weak non-rigidity, that is $ \boldsymbol{\Pi}_{\Delta_{0,0}}(\boldsymbol{\Gamma}^r)=q\ex^{-|\Delta_{0,0}|}|\Delta_{0,0}|z>0$.
%
Using \cite[Theorem~3.3]{DDG12} and \cite[Corollary~3.4]{DDG12} we obtain all the statements in Proposition~\ref{Propexist}. \qed

\subsection{Breaking of the symmetry of the mark distribution}\label{finish}
In this section we complete the proof of Theorem~\ref{THM-mainedge}  by analysing the Gibbs distributions $ \gamma_{\L,\bo} $ in the limit $ \L\uparrow\R^2 $. We pick an admissible boundary condition $ \om\in\O^*_{\L_n} $,  and we  let $ \bo=(\om\setminus\L_n,\emptyset,\ldots,\emptyset) $ be the admissible monochromatic boundary condition such that $ \bo\in\bO_{\L_n}^* $. We write $ \gamma_n $ for $ \gamma_{\L_n,\bo} $ and  let $ P_n $ be the probability measure on $ \bO $ relative to which the marked configurations in distinct parallelotopes $ \L_n+(2n+1)M(k,l), (k,l)\in\Z^2 $, are independent with identical distribution $ \gamma_n $. As we are dealing with a cell structure for the partition of $ \R^2 $,  we confine ourself first to lattice shifts when we employ spatial averaging. Thus,
$$
\overline{P}_n=\frac{1}{2n+1}\sum_{(k,l)\in\{-n,\ldots,n\}^2}P_n\circ\theta_{M(k,l)}^{-1}.
$$
By the periodicity of $ P_n $ the measure $\overline{P}_n $ is $\Z^2$-shift-invariant. The proof in \cite[Chapter 5]{DDG12} shows that $ (\overline{P}_n)_{n\ge 1} $ has a subsequence which converges with respect to the topology of local convergence to some $ \widehat{P} \in\Mcal_1(\bO) $. As outlined in \cite{DDG12} it is difficult to show that $ \widehat{P} $ is concentrated on admissible configurations. As $ \widehat{P} $ is non-degenerate the proof in \cite[Chapter 5]{DDG12} shows that $ P=\widehat{P}(\cdot|\{\emptyset\}^{\rm c}) $ is a Gibbs measure with $ P(\{\emptyset\})=0 $.  In order to obtain an $ \R^2$-shift-invariant Gibbs measure one needs to apply another averaging,
$$
P^{\ssup{1}}=\int_{\Delta_{0,0}}\,P\circ\theta_{Mx}^{-1}\,\d x.
$$
Applying Propositions~\ref{Prop-sym} and \ref{Prop-per}, we see that for $ \Delta=\Delta_{0,0} $,
$$
\begin{aligned}
\int\,(qN_{\Delta,1}-N_{\Delta})\,\d \overline{P}_n& \ge \frac{(q-1)}{2n+1}\sum_{(k,l)\in\{-n,\ldots,n\}^2}\int\,N_{\Delta_{k,l}\leftrightarrow\infty}\,\d C_{\L_n,\om}\\
&\ge (q-1)\eps.
\end{aligned}
$$ Thus
$$
\int\,(qN_{\Delta,1}-N_\Delta)\,\d P^{\ssup{1}}>0,
$$ and  we observe the following break of symmetry in the expected density of particles of type $1 $ and of any other  type, that is, 
$$
\rho_1(P^{\ssup{1}})>\rho_2(P^{\ssup{1}})=\cdots=\rho_q(P^{\ssup{1}}),
$$ where $ \rho_s(P^{\ssup{1}})=1/|\Delta|\E_{P^{\ssup{1}}}[N_{\Delta,s}], s\in M_q $. We conclude with our statement as in \cite{GH96} by showing that the matrix 
$$
\big(\rho_s(P^{\ssup{t}})\big)_{s,t\in M_q}
$$ is regular, where $ P^{\ssup{t}} $ is obtained from $ P^{\ssup{1}} $ by swapping the role of $1$ and $ t $.
\section*{Appendix}
\begin{appendices}

 \section{Pseudo-periodic configurations}\label{pseudoperiodic}
  We define pseudo-periodic configurations as in \cite{DDG12}. We obtain a partition of $ \R^2 $ which is adapted to the Delaunay tessellation. Pick a length scale $ \ell>0 $ and consider the matrix $$ M=\left(\begin{matrix}M_1 & M_2 \end{matrix}\right) =\left(\begin{matrix} \ell & \ell/2\\0 &\sqrt{3}/2 \ell\end{matrix}\right).$$
 Note that $ |M_i|=\ell, i=1,2 $, and $ \angle(M_1,M_2)=\pi/3 $. For each $ (k,l)\in\Z^2 $ we define the cell
 \begin{equation}\label{cell}
 \Delta_{k,l}=\{Mx\in\R^2\colon x-(k,l)\in[-1/2,1/2)^2\}
 \end{equation} with area $ |\Delta_{k,l}|=\frac{\sqrt{3}}{2}\ell^2 $. These cells constitute a periodic partition of $ \R^2 $ into parallelotopes. Let $ B $ be a measurable set of $\O_{\Delta_{0,0}}\setminus\{\emptyset\} $ and
 \begin{equation}\label{ppc}
\Gamma=\{\om\in\O\colon\theta_{Mz}(\om_{\Delta_{k,l}})\in B\mbox{ for all } (k,l)\in\Z^2\} 
 \end{equation} the set of all configurations whose restriction to an arbitrary cell, when shifted back to $ \Delta_{0,0} $, belongs to $ B $. Elements of $ \Gamma $ are called \textbf{pseudo-periodic} configurations. We define marked pseudo-periodic configurations in an analogous way.

\section{Topology of local convergence}\label{AppC}
We write $ \Mcal_1^{\varTheta}(\bO) $ (resp. $ \Mcal_1^{\varTheta}(\O) $) for the set of all shift-invariant probability measures on $ (\bO,\boldsymbol{\Fcal}) $ (resp. $ (\O,\Fcal) $). A measurable function $ f\colon\bO\to\R $ is called local and tame if
$$
f(\bo)=f(\bo_\L)\quad\mbox{ and }\quad |f(\bo)|\le a N_\L(\bo)+b
$$ for all $ \bo\in\bO $ and some $ \L\Subset\R^2 $ and suitable constants $ a,b\ge 0 $. Let $ \Lscr $ be the set of all local and tame functions. The topology of local convergence, or $\Lscr$-topology, on $ \Mcal_1^{\varTheta}(\bO) $ is then defined as the weak$*$ topology induced by $ \Lscr $, i.e., as the smallest topology for which the mappings $ P\mapsto \int f\d P $ with $ f\in\Lscr $ are continuous.

\section{Mixed site-bond percolation}\label{AppD}
Given a graph $ G=(V,E) $, let $ \P_p $ be the probability measure on configurations of open and closed sites of $G$. Each site of $G$ is open with probability $p $ and closed with probability $ 1-p $. Similarly, let $\widetilde{\P}_p $ be the probability measure on configurations of open and closed edges of $G$. Each edge of $G$ is open with probability $p$ and closed with probability $ 1-p$. For $ x_0\in V $ and a subset of vertices $X\subset V $, let
$$
\begin{aligned}
\bs(p,x_0,X,G)&=\P_p(\exists\,\mbox{ a path } x_0=v_0,e_1,v_1,\ldots, e_n,v_n\,\mbox{ with }\,  v_n\in X\,\mbox{  and all vertices are open}),\\
\bbeta(p,x_0,X,G)&=\widetilde{\P}_p(\exists\,\mbox{ a path } x_0=v_0,e_1,v_1,\ldots, e_n,v_n\,\mbox{ with }\,  v_n\in X\,\mbox{  and all edges are open}).\end{aligned}
$$
It is known since \cite{Kesten} that site percolation implies bond percolation, that is, for any $ 0\le p\le 1$,
\begin{equation}\label{mixed}
\bs(p,x_0,X,G)\le \bbeta(p,x_0,X,G).
\end{equation}

In mixed site-bond percolation, both edges and vertices may be open or closed, possibly with different probabilities. Each edge or bond is open independently of anything else with probability $ p^\prime $ and each site is open independently of anything else with probability $ p $. The edges and sites that are not open, along with the edges to or from these sites, are closed. We shall consider paths of open sites and open edges. For $x_0\in V$ and a subset of vertices $X\subset V $, let
$$
\begin{aligned}
\bgamma(p,p^\prime,x_0,X,G)&=\P_{pp^\prime}(\exists\,\mbox{ a path } x_0=v_0,e_1,v_1,\ldots, e_n,v_n\,\mbox{ with }\,  v_n\in X\,\\ & \qquad
\mbox{  and all vertices and all edges are open}).
\end{aligned}
$$ 
Let $G^\prime $ be the reduced graph where each edge and site of  $ G $ is removed independently with probability $1-p^\prime $ and $1-p $ respectively. By taking  the expectation on both sides of inequality \eqref{mixed}, on $ G^\prime $, with respect to $ \P_\delta $ and $ \widetilde{P}_\delta $, we arrive at the mixed site-bond percolation result of Hammersley, a generalistion of the work of McDiarmid, see  \cite{Ham}. That is, for $\delta,p,p^\prime\in [0,1] $ one gets that
\begin{equation}\label{eqHam}
\bgamma(\delta p,p^\prime,x_0,X,G) \le   \bgamma(p,\delta p^\prime,x_0,X,G).
\end{equation}

By setting $ \delta=p $ and $ p^\prime =1 $ in \eqref{eqHam}, and noticing that $ \bgamma(p^2,1,x_0,X,G)=\bs(p^2,X_0,X,G) $, we arrive at
\begin{equation}
\bs(p^2,x_0,X,G)\le \bgamma(p,p,x_0,X,G),
\end{equation}
and hence
\begin{equation}\label{ineqmixed}
\theta_{\ssup{\rm mixed}}(p,p)\ge \theta_{\ssup{\rm site}}(p^2),
\end{equation}
where $ \theta_{\ssup{\rm mixed}}(p,p^\prime) $ is the mixed site-bond percolation probability with parameters $p $ and $ p^\prime $, and $ \theta_{\ssup{\rm site}}(p) $ is the site percolation probability with parameter $ p $.

\section{Geometrical Lemmas}\label{appE}
\begin{lemma}\label{L:E} Let $ \Gamma=(V,E) $ be a spoked chain with $ V=\{x_1,\ldots, x_n\} $ and $ \hat{x}_1<\cdots<\hat{x}_n $. For $ 1<k\le n $, let $ x^*_k $ and $ x^*_{k+1} $ be the centres of the circumscribing circles of the triangles $ \tau(x_0,x_{k-1},x_k) $ and $ \tau(x_0,x_k,x_{k+1}) $ respectively. Then $ \hat{x}^*_{k+1}\ge \hat{x}^*_k $.
\end{lemma}
\begin{proofsect}{Proof}
The points $ x^*_{k+1} $ and $ x^*_k $ both lie on the bisector of the line segment $ \overline{x_0x_k} $. Suppose $ \hat{x}^*_k>\hat{x}_k $, then the radius of circumcircle $\partial B(\tau(x_0,x_k,x_{k+1})) $ is greater than the radius of $ B(\tau(x_0,x_{k-1},x_k)) $ and hence $ \hat{x}^*_{k+1}\ge \hat{x}^*_k $. Now suppose $ \hat{x}^*_k\le \hat{x}_k $. If $ \hat{x}^*_{k+1}<\hat{x}^*_k $, then $ x_{k+1} $ lies in the interior of the circumcircle $ \partial B(\tau(x_0,x_{k-1},x_k)) $ which contradicts properties of the Delaunay tessellation. Therefore, $ \hat{x}^*_{k+1}\ge \hat{x}^*_k $.
\qed
\end{proofsect}
Let $ a\in\R^2 $ be the pole in a polar coordinate system where $ \hat{x} $ denotes the angular coordinate of $ x\in\R^2 $. For $x,y\in\R^2 $ with $ \hat{x}<\hat{y} $, let $ \partial B(\tau(a,x,y)) $ be the unique circle that circumscribes all three vertices (Delaunay tessellation). Let $  C^a_{xy} $ be the arc opposite the vertex $ a $. For any arc $C$, let $L(C) $ denote its length. 
\begin{theorem}\label{T:E}
Suppose $a \in\R^2 $ is the pole. Let $ b,c\in\R^2 $ with $ 0<\hat{b}<\hat{c}<\pi $. Let $U$ be the convex hull of $ C^a_{bc} $. Then, for all $z\in U $,
\begin{equation}
L(C^a_{bz})+L(C^a_{zc})\le L(C^a_{bc}).
\end{equation}
\end{theorem}
\begin{proofsect}{Proof}
Let $ r>0 $ denote the radius of the circumcircle $ \partial B(\tau(a,b,c)) $ and  define for $ z\in U $, 
$$
\begin{aligned}
M&:=|b-c|;\; h_1:=|b-z|; \;h_2:=|z-c|;\;t:=|z-a|;\; s_1:=|b-a|;\; s_2:=|c-a|,\\
\theta_1&:=\hat{z}-\hat{b};\;\theta_2:=\hat{c}-\hat{z};\;\theta:=\theta_1+\theta_2.
\end{aligned}
$$
Then, $ L(C^a_{bz})=2\theta\radius(B(\tau(a,b,z))) $ with $ \radius(B(\tau(a,b,z)))=h_1/2\sin(\theta_1) $. Thus the following holds:
$$
L(C^a_{bz})=h_1\frac{\theta_1}{\sin(\theta_1)}, \; L(C^a_{zc})=h_2\frac{\theta_2}{\sin(\theta_2)},\; L(C^a_{bc})=M\frac{\theta}{\sin(\theta)}.
$$
The strategy of the proof is to first show that $ L(C^a_{bz})+L(C^a_{zc})=L(C^a_{bc}) $ for $ z\in C^a_{bc} $ and $ L(C^a_{bz})+L(C^a_{zc})\le L(C^a_{bc}) $ for $ z\in\overline{bc} $.
We then define $  L(C^a_{bz})+L(C^a_{zc}) $ as a function of $ \theta_1,s_1,t $ and $ r $, and show that it is convex with respect to $t $. Noting that $z\in U $ is uniquely determined by $ t $ and $ \theta_1 $, we conclude with the result for all $ z\in U $.\\[1ex]
Let $ z\in C^a_{bc} $. Then $ B(\tau(a,b,c))=B(\tau(a,b,z))=B(\tau(a,z,c)) $. Therefore, $ C^1_{bz}\cup C^a_{zc}=C^a_{bc} $ and thus \begin{equation}\label{E2} L(C^a_{bz})+L(C^a_{zc})=L(C^a_{bc}).\end{equation} 
Now let $ z\in\partial U\cap\overline{bc} $. Then, $ h_1+h_2=M $ and
\begin{equation}\label{E3}
\begin{aligned}
L(C^a_{bz})+L(C^a_{zc})&=h_1\frac{\theta_1}{\sin(\theta_1)}+h_2\frac{\theta_2}{\sin(\theta_2)}=h_1\frac{\theta_1}{\sin(\theta_1)}+(M-h_1)\frac{\theta_2}{\sin(\theta_2)}\\
&\le h_1\frac{\theta}{\sin(\theta)}+(M-h_1)\frac{\theta}{\sin(\theta)}=M\frac{\theta}{\sin(\theta)}=L(C^a_{bc}),
\end{aligned}
\end{equation}
where the inequality holds because $ \theta\ge\max\{\theta_1,\theta_2\} >0$ and $ g(x):=\frac{x}{\sin(x)} $ is an increasing function on the interval $ [0,\pi]$. To write $ L(C^a_{bz})+L(C^a_{zc}) $ as a function of $ \theta_1,s_1 $ and $t$, note that by the cosine rule of triangles,
$$
h_1^2=t^2-2s_1t\cos(\theta_1)+s_1^2 \;\mbox{ and }\; h_2^2=t^2-2s_2t\cos(\theta_2)+s_2^2,
$$
and thus
$$
L(C^a_{bz})=(t^2-2s_1t\cos(\theta_1)+s_1^2)^{1/2}\frac{\theta_1}{\sin(\theta_1)}=:f_1(\theta_1,s_1,t).
$$
Furthermore, $ s_2 $ is a function of $s_1 $ and $\theta_1 $ since
$$
M^2=s_1^2+s_2^2-2s_1s_2\cos\big(\sin^{-1}\big(\frac{M}{2r}\big)\big),
$$
and $ \theta_2 $ is a function of $ \theta_1 $, $$ \theta_2=\theta-\theta_1=\sin^{-1}\big(\frac{M}{2r}\big)-\theta_1. 
$$
We obtain from these relations the expression
\begin{equation}
L(C^a_{zc})=\big(t^2-2s_2(s_1,\theta_1)t\cos(\theta_2)+s_2(s_1,\theta_1)^2\big)^{1/2}\frac{\theta_2(\theta_1)}{\sin(\theta_2(\theta_1))}=:f_2(\theta_1,s_1,t).
\end{equation}
We will show that $ f(\theta_1,s_1,t)=f_1(\theta_1,s_1,t)+f_2(\theta_1,s_1,t) $ is convex with respect to $t$. We obtain
$$
\begin{aligned}
&\frac{\d^2}{\d t^2}f_1(\theta_1,s_1,t)=\\
&=\frac{\theta_1}{\sin(\theta_1)}\Big(\frac{(t^2-2s_1t\cos(\theta_1)+s_1^2)^3/2}{(t^2-2s_1t\cos(\theta_1)+s_1^2)^2}-\frac{(t^2-2s_1t\cos(\theta_1)+s_1^2)^{1/2}(t-s_1\cos(\theta_1))^2}{(t^2-2s_1t\cos(\theta_1)+s_1^2)^2}\Big).\end{aligned}
$$
The function $ \frac{x}{\sin(x)} $ is positive for $ 0\le x\le \pi $. The denominator in the bracket is just $ h_1^2 $ and thus positive. The numerator in the bracket reads as
$$
\begin{aligned}
& t^2-2s_1t\cos(\theta_1)+s_1^2)^3/2-(t^2-2s_1t\cos(\theta_1)+s_1^2)^{1/2}(t-s_1\cos(\theta_1))^2\\
&=(t^2-2s_1t\cos(\theta_1)-s_1^2)^{1/2}(t^2-2s_1t\cos(\theta_1)+s_1^2-(t-s_1\cos(\theta_1))^2)\\
&=(t^2-2s_1t\cos(\theta_1)-s_1^2)^{1/2} s_1^2(1-\cos^2(\theta_1))\ge 0,
\end{aligned}
$$
since $(t^2-2s_1t\cos(\theta_1)-s_1^2)^{1/2}=h_1\ge 0$. Therefore, the function $f_1 $ is convex with respect to $ t$. Similarly, show that $ f_2 $ is convex with respect to $t $ to see that the function $f$ is convex with respect to $ t $. 
Pick $ 0\le \theta_1\le\sin^{-1}\big(\frac{M}{2r}\big) $. There exist $ 0<t_{\rm min}(\theta_1)<t_{\rm max}(\theta_1)<2r $ such that $ t_{\rm min}(\theta_1)\le |z|\le t_{\rm max}(\theta_1) $ for all $z\in U $ with $\hat{z}-\hat{b}=\theta_1 $. We have shown (see \eqref{E2} and \eqref{E3}) that
$$
f(\theta_1,s_1,t_{\rm min}(\theta_1))\le L(C^a_{bc})\;\mbox{ and }\; f(\theta_1,s_1,t_{\rm max}(\theta_1))=L(C^a_{bc}).
$$
Therefore, by the convexity of $f$, for all $ t\in [t_{\rm min}(\theta_1),t_{\rm max}(\theta_1)]$,
$$
\begin{aligned}
f(\theta_1,s_1,t)&\le\frac{t-t_{\rm min}(\theta_1)}{t_{\rm max}(\theta_1)-t_{\rm min}(\theta_1)}f(\theta_1,s_1,t_{\rm min}(\theta_1))+\frac{t_{\rm max}(\theta_1)-t}{t_{\rm max}(\theta_1)-t_{\rm min}(\theta_1)}f(\theta_1,s_1,t_{\rm max}(\theta_1))\\&\le f(\theta_1,s_1,t_{\rm max}(\theta_1))=L(C^a_{bc}).
\end{aligned}
$$
\qed
\end{proofsect}

\end{appendices}
\section*{Acknowledgments}
S.A. thanks  D. Dereudre for helpful discussions and exchange of ideas visiting Lille.

\end{document}